\documentclass{siamltex1213}
\usepackage[paperwidth=7in, paperheight=10in, margin=.875in]{geometry}
\usepackage{amsfonts,amssymb}
\usepackage{amsmath}
\usepackage{graphicx}
\usepackage{enumerate}
\usepackage{pgf,tikz}
\usetikzlibrary{arrows, automata}
\usepackage{mathrsfs}
\usepackage{lettrine}
\usepackage{supertabular}
\usepackage{psfrag}
\usepackage{curves,eepic}
\usepackage{xspace}
\usepackage{ifthen}
\usepackage{rotating}
\usepackage{fancyvrb}
\usepackage{epic}
\usepackage{eepic}
\usepackage{enumerate}
\usepackage{multicol}
\usepackage{amsmath,amssymb,upref,amscd}
\usepackage{epsfig}
\usepackage{subcaption}
\usepackage[T1]{fontenc}
\usepackage{textcomp}
\usepackage{lmodern}
\usepackage[sc]{mathpazo}
\usepackage{verbatim}
\usepackage{algorithm}
\usepackage{algorithmic}
\usepackage{hyperref}
\usepackage{url}
\usepackage{color}
\usepackage{colordvi}
\usepackage{tabularx}
\usepackage{overpic}
\usepackage{tikz}
\usepackage{pgfplots}
\usetikzlibrary{spy,calc}

\newif\ifblackandwhitecycle
\gdef\patternnumber{0}

\pgfkeys{/tikz/.cd,
    zoombox paths/.style={
        draw=orange,
        very thick
    },
    black and white/.is choice,
    black and white/.default=static,
    black and white/static/.style={ 
        draw=white,   
        zoombox paths/.append style={
            draw=white,
            postaction={
                draw=black,
                loosely dashed
            }
        }
    },
    black and white/static/.code={
        \gdef\patternnumber{1}
    },
    black and white/cycle/.code={
        \blackandwhitecycletrue
        \gdef\patternnumber{1}
    },
    black and white pattern/.is choice,
    black and white pattern/0/.style={},
    black and white pattern/1/.style={    
            draw=white,
            postaction={
                draw=black,
                dash pattern=on 2pt off 2pt
            }
    },
    black and white pattern/2/.style={    
            draw=white,
            postaction={
                draw=black,
                dash pattern=on 4pt off 4pt
            }
    },
    black and white pattern/3/.style={    
            draw=white,
            postaction={
                draw=black,
                dash pattern=on 4pt off 4pt on 1pt off 4pt
            }
    },
    black and white pattern/4/.style={    
            draw=white,
            postaction={
                draw=black,
                dash pattern=on 4pt off 2pt on 2 pt off 2pt on 2 pt off 2pt
            }
    },
    zoomboxarray inner gap/.initial=5pt,
    zoomboxarray columns/.initial=2,
    zoomboxarray rows/.initial=2,
    subfigurename/.initial={},
    figurename/.initial={zoombox},
    zoomboxarray/.style={
        execute at begin picture={
            \begin{scope}[
                spy using outlines={%
                    zoombox paths,
                    width=\imagewidth / \pgfkeysvalueof{/tikz/zoomboxarray columns} - (\pgfkeysvalueof{/tikz/zoomboxarray columns} - 1) / \pgfkeysvalueof{/tikz/zoomboxarray columns} * \pgfkeysvalueof{/tikz/zoomboxarray inner gap} -\pgflinewidth,
                    height=\imageheight / \pgfkeysvalueof{/tikz/zoomboxarray rows} - (\pgfkeysvalueof{/tikz/zoomboxarray rows} - 1) / \pgfkeysvalueof{/tikz/zoomboxarray rows} * \pgfkeysvalueof{/tikz/zoomboxarray inner gap}-\pgflinewidth,
                    magnification=3,
                    every spy on node/.style={
                        zoombox paths
                    },
                    every spy in node/.style={
                        zoombox paths
                    }
                }
            ]
        },
        execute at end picture={
            \end{scope}
     \gdef\patternnumber{0}
        },
        spymargin/.initial=0.5em,
        zoomboxes xshift/.initial=1,
        zoomboxes right/.code=\pgfkeys{/tikz/zoomboxes xshift=1},
        zoomboxes left/.code=\pgfkeys{/tikz/zoomboxes xshift=-1},
        zoomboxes yshift/.initial=0,
        zoomboxes above/.code={
            \pgfkeys{/tikz/zoomboxes yshift=1},
            \pgfkeys{/tikz/zoomboxes xshift=0}
        },
        zoomboxes below/.code={
            \pgfkeys{/tikz/zoomboxes yshift=-1},
            \pgfkeys{/tikz/zoomboxes xshift=0}
        },
        caption margin/.initial=4ex,
    },
    adjust caption spacing/.code={},
    image container/.style={
        inner sep=0pt,
        at=(image.north),
        anchor=north,
        adjust caption spacing
    },
    zoomboxes container/.style={
        inner sep=0pt,
        at=(image.north),
        anchor=north,
        name=zoomboxes container,
        xshift=\pgfkeysvalueof{/tikz/zoomboxes xshift}*(\imagewidth+\pgfkeysvalueof{/tikz/spymargin}),
        yshift=\pgfkeysvalueof{/tikz/zoomboxes yshift}*(\imageheight+\pgfkeysvalueof{/tikz/spymargin}+\pgfkeysvalueof{/tikz/caption margin}),
        adjust caption spacing
    },
    calculate dimensions/.code={
        \pgfpointdiff{\pgfpointanchor{image}{south west} }{\pgfpointanchor{image}{north east} }
        \pgfgetlastxy{\imagewidth}{\imageheight}
        \global\let\imagewidth=\imagewidth
        \global\let\imageheight=\imageheight
        \gdef\columncount{1}
        \gdef\rowcount{1}
        
    },
    image node/.style={
        inner sep=0pt,
        name=image,
        anchor=south west,
        append after command={
            [calculate dimensions]
            node [image container,subfigurename=\pgfkeysvalueof{/tikz/figurename}-image] {\phantomimage}
            node [zoomboxes container,subfigurename=\pgfkeysvalueof{/tikz/figurename}-zoom] {\phantomimage}
        }
    },
    color code/.style={
        zoombox paths/.append style={draw=#1}
    },
    connect zoomboxes/.style={
    spy connection path={\draw[draw=none,zoombox paths] (tikzspyonnode) -- (tikzspyinnode);}
    },
    help grid code/.code={
        \begin{scope}[
                x={(image.south east)},
                y={(image.north west)},
                font=\footnotesize,
                help lines,
                overlay
            ]
            \foreach \x in {0,1,...,9} { 
                \draw(\x/10,0) -- (\x/10,1);
                \node [anchor=north] at (\x/10,0) {0.\x};
            }
            \foreach \y in {0,1,...,9} {
                \draw(0,\y/10) -- (1,\y/10);                        \node [anchor=east] at (0,\y/10) {0.\y};
            }
        \end{scope}    
    },
    help grid/.style={
        append after command={
            [help grid code]
        }
    },
}

\newcommand\phantomimage{%
    \phantom{%
        \rule{\imagewidth}{\imageheight}%
    }%
}
\newcommand\zoombox[2][]{
    \begin{scope}[zoombox paths]
        \pgfmathsetmacro\xpos{
            (\columncount-1)*(\imagewidth / \pgfkeysvalueof{/tikz/zoomboxarray columns} + \pgfkeysvalueof{/tikz/zoomboxarray inner gap} / \pgfkeysvalueof{/tikz/zoomboxarray columns} ) + \pgflinewidth
        }
        \pgfmathsetmacro\ypos{
            (\rowcount-1)*( \imageheight / \pgfkeysvalueof{/tikz/zoomboxarray rows} + \pgfkeysvalueof{/tikz/zoomboxarray inner gap} / \pgfkeysvalueof{/tikz/zoomboxarray rows} ) + 0.5*\pgflinewidth
        }
        \edef\dospy{\noexpand\spy [
            #1,
            zoombox paths/.append style={
                black and white pattern=\patternnumber
            },
            every spy on node/.append style={#1},
            x=\imagewidth,
            y=\imageheight
        ] on (#2) in node [anchor=north west] at ($(zoomboxes container.north west)+(\xpos pt,-\ypos pt)$);}
        \dospy
        \pgfmathtruncatemacro\pgfmathresult{ifthenelse(\columncount==\pgfkeysvalueof{/tikz/zoomboxarray columns},\rowcount+1,\rowcount)}
        \global\let\rowcount=\pgfmathresult
        \pgfmathtruncatemacro\pgfmathresult{ifthenelse(\columncount==\pgfkeysvalueof{/tikz/zoomboxarray columns},1,\columncount+1)}
        \global\let\columncount=\pgfmathresult
        \ifblackandwhitecycle
            \pgfmathtruncatemacro{\newpatternnumber}{\patternnumber+1}
            \global\edef\patternnumber{\newpatternnumber}
        \fi
    \end{scope}
}

\renewcommand{\theequation}{\arabic{section}.\arabic{equation}}

\newtheorem{problem}{Problem}
\newtheorem{remark}[theorem]{Remark}
\DeclareMathOperator{\Ann}{ann}
\DeclareMathOperator{\esssup}{ess\ sup}

\DeclareMathOperator{\diver}{div}

\DeclareMathOperator*{\argmin}{arg\,min}
   \allowdisplaybreaks
\begin{document}
 \title{Iterative TV Minimization on the Graph}



\author{Japhet Niyobuhungiro\thanks{Department of Mathematics, School of Science, College of Science and Technology, University of Rwanda, P.O. Box 3900 Kigali, Rwanda, (japhetniyo@gmail.com, jniyobuhungiro@ur.ac.rw).}
\and Eric Setterqvist\thanks{Faculty of Mathematics, University of Vienna, Oskar-Morgenstern-Platz 1, A-1090 Vienna, Austria, (eric.setterqvist@univie.ac.at).}
\and Freddie \r{A}str\"om\thanks{Heidelberg University, IWR/R.B108, Berliner Str. 43, 69120 Heidelberg, Germany, (freddie.astroem@iwr.uni-heidelberg.de, freddie.astrom@gmail.com).}
\and George Baravdish\thanks{Department of Science and Technology, Link{\"o}ping University, SE-601 74 Norrk\"oping, Sweden, (george.baravdish@liu.se).}}

\pagestyle{myheadings} \markboth{Iterative TV Minimization on the Graph}{J. Niyobuhungiro, E. Setterqvist, F. \r{A}str\"om \& G. Baravdish} \maketitle
\begin{abstract}
We define the space of functions of bounded variation ($BV$) on the graph. Using the notion of divergence of flows on graphs, we show that the unit ball of the dual space to $BV$ in the graph setting can be described as the image of the unit ball of the space $\ell^{\infty}$ by the divergence operator. Based on this result, we propose a new iterative algorithm to find the exact minimizer for the total variation (TV) denoising problem on the graph. The proposed algorithm is provable convergent and its performance on image denoising examples is compared with the Split Bregman and Primal-Dual algorithms as benchmarks for iterative methods and with BM3D as a benchmark for other state-of-the-art denoising methods. The experimental results show highly competitive empirical convergence rate and visual quality for the proposed algorithm.
\end{abstract}
\section{Introduction}\label{sec:Intro}
\subsection{Background}
Removing or reducing the noise from obtained and observed images is a fundamental image processing problem known as \textit{denoising} appearing in many application areas. The image noise $\eta$ considered here is additive which means that the observed image data $u_{0}$ is related to the underlying true image $u$ according to the linear model
\begin{align} \label{eq:linmodel}
  u_{0} = u + \eta.
\end{align}
The noise component $\eta$ is further assumed to be normally, independent and identically distributed. We study in this work the total variation (abbreviated as TV) denoising problem. In this problem a noisy image $u_{0}\in L^{2}\left(\Omega\right)$, where the open set $\Omega\subset\mathbb{R}^{2}$ is the image domain, is observed and the denoised image approximating the original image is then defined as the solution $u_{opt}$ of the optimization problem
\begin{align}
\label{TheMinProblem}
	\inf_{u\in BV(\Omega)} \left( \frac{1}{2}\left\|u_{0}-u\right\|^{2}_{L^{2}(\Omega)}+t\left\|u\right\|_{BV(\Omega)} \right),
\end{align}
where $t>0$ is called the regularization parameter and $BV\left(\Omega\right)$ is the space of functions of bounded variation. The TV denoising model \eqref{TheMinProblem} was introduced in 1992 by L. I. Rudin, S. Osher and E. Fatemi \cite{ROFPaper} and is now also widely known in the image processing community as the ROF model.
The space $BV(\Omega)$ is defined as follows
\begin{definition}
\begin{align*}
BV(\Omega)=\left\{u\in L^{1}(\Omega):\left\|u\right\|_{BV}<\infty\right\}
\end{align*}
\end{definition}
where the bounded variation or total variation seminorm of $u$ is given by
\begin{definition}
 \begin{align*}
	 \left\|u\right\|_{BV(\Omega)} = \int_{\Omega}\left|Du\right|=\sup \int_{\Omega}u(x)\diver\ \vec{g}(x)dx,
 \end{align*}
where the supremum is taken over all $\vec{g}\in C^{1}_{c}\left(\Omega,\mathbb{R}^{2}\right)$ such that $\sup_{x\in\Omega}\sqrt{g^{2}_{1}(x)+g^{2}_{2}(x)}\leq 1$.
\end{definition}
Note that if $u$ is a differentiable function then $\left\|u\right\|_{BV(\Omega)} = \int_{\Omega}\left|\nabla u(x)\right|dx$.
An important feature of the $BV$ term in the minimization problem \eqref{TheMinProblem} is that it discourages the solution from having oscillations and at the same time allowing it to have discontinuities. 
\par Since its appearance in 1992, the ROF model has received a large amount of popularity for its effeciency in denoising images without smoothing out the boundaries, and it has also been applied to a multitude of other imaging problems (see for example the book \cite{TonyFChan}). We choose next to highlight a few selected works from the vast literature on the ROF model and TV minimization which are related to our approach.
\par An early work on total variation minimization based on dual formulation is \cite{OFChambolle12}. In 2004, A. Chambolle provided an iterative algorithm related to \cite{OFChambolle12} and proved its convergence, see \cite{Chambolle2004}. We remark that the works \cite{AujolChamGLBF2005, ChambolleDarbon2006} also proposed efficient projection algorithms for total variation minimization. The papers \cite{JEROMEglesiSiam} and \cite{JFAGuyGTFCSO2006} adapted Chambolle's algorithm from \cite{Chambolle2004} to handle linear operators in the ROF model, such as convolution operators representing blurring.
\par After the appearance of \cite{Chambolle2004}, several iterative algorithms have been developed which can be used to solve TV minimization problems. Bregman iteration was shown in \cite{TomGoldsteinStanleyOsher} to be a efficient and fast way to solve TV problems among other $L^{1}$-regularized optimization problems. In particular, a \emph{split Bregman} method was proposed in \cite{TomGoldsteinStanleyOsher} and subsequently used to compute the ROF minimizer. The Primal-Dual algorithm proposed in \cite{AChambollePock} is another general purpose iterative algorithm which can be efficiently applied to solve TV minimization problems. The fast iterative shrinkage-thresholding algorithm (FISTA) for linear inverse problems, see \cite{BeckTeboulleFista1} and \cite{BeckTeboulleFista2}, is also known to be able to solve TV minimization problems efficiently.
\par For anisotropic total variation minimization of quantized images, i.e. the pixels of the image take values in a prescribed finite set because the observed image is decomposed into a prescribed number of level sets, graph cut algorithms have been developed that exactly compute the minimizer up to machine precision. Foundational works in this direction are the algorithms of Chambolle \cite{AChambBinInt}, Darbon and Sigelle \cite{DarbonSigelle2006} and Goldfarb and Yin \cite{DWYin2009}. These algorithms are not iterative and in terms of speed, they are very fast. 
\par Based on the fact that an image has a locally sparse representation in transform domain and that this sparsity is enhanced by grouping similar 2D image patches into 3D groups, a paper on Collaborative filtering or BM3D grouping and filtering procedure was written \cite{BM3D1} and later analysed and implemented in \cite{BM3D2}. Though this method is not directly designed to solve the TV minimization problem, it is one of the state-of-the-art filtering methods applicable to the denoising problem.
\par Image decomposition models into a piecewise-smooth and oscillating components that usually researchers refer to as cartoon and textures (or textures $+$ noise) respectively, have received great interest in the image processing community. For example $u_{0}\in L^{2}\left(\Omega\right)$ is decomposable as
\begin{align*}
  u_{0} = u_{opt} + \left( u_{0}-u_{opt} \right).
\end{align*}
This is the decomposition of $u_{0}$ into the piecewise-smooth component $u_{opt}\in BV\left(\Omega\right)$ satisfying \eqref{TheMinProblem} and the component $\left(u_{0}-u_{opt}\right)\in L^{2}\left(\Omega\right)$ which contains textures and noise. The original theoretical model for such an image decomposition was introduced in 2001 by Y. Meyer in \cite{YvesMeyer} by using the total variation to model the piecewise-smooth component and an appropriate dual space $G$ which is the Banach space composed of the distributions $f=\partial_{1}g_{1}+\partial_{2}g_{2}=\diver \vec{g}$, where $g_{1}$ and $g_{2}$ are in $L^{\infty}(\Omega)$ and $\left\|f\right\|_{G} = \inf \left\|\vec{g}\right\|_{L^{\infty}(\Omega;\mathbb{R}^{2})}$ where the infimum is taken over all $\vec{g}$ such that $f = \diver \vec{g}$ and $\left\|\vec{g}\right\|_{L^{\infty}(\Omega;\mathbb{R}^{2})}=\esssup_{x\in\Omega}\sqrt{\left|g_{1}(x)\right|^{2}+\left|g_{2}(x)\right|^{2}}$, to model the oscillating component. Some of the works proposed in the literature for numerically solving Meyer's model or its variants include for instance \cite{AujolChambolle2005} that proposed to split the image into three components, a geometrical component modeled by the total variation, a texture component modeled by a negative Sobolev norm and a noise component modeled by a negative Besov norm.  Furthermore, \cite{JGillesSOsher} designed an algorithm by using split Bregman iterations and the duality used by Chambolle to find the minimizer of a functional based on Meyer's $G$-norm. Other works based on the $G$-norm include for example \cite{DMStrongAujolTFC2006} and \cite{OsherScherzerGNorm}.
\subsection{Summary of main contributions and motivation}
We present an iterative method for solving the discrete analogue of the TV minimization problem \eqref{TheMinProblem} on finite graphs. The algorithmic representation of the method is given in Algorithm \ref{ALGO} and is proved to converge to the exact minimizer. Further, the algorithm can be run on a parallel computer architecture and is thereby suitable to handle large graphs and data sets. The proof of the convergence result Theorem \ref{PropoTconv} is based on duality principles from convex analysis and Theorem \ref{BVstarBall} which characterizes, in the graph setting, the unit ball of the dual space to $BV$ as the image of the unit ball of the space $\ell^{\infty}$ by the divergence operator. We note that the strength of a graph representation is when considering non-Euclidean metric spaces via manifold representations, for example when the image is a map in spherical geometry, which could be the case in many applications. Our approach also illustrates the properties of the optimal decomposition of the image data into a piecewise-smooth image component and a noise component and gives its geometrical interpretation. Experimental results confirm that our method is a highly competitive TV denoising algorithm in terms of both convergence rate and visual quality.
\subsection{Overview}
This paper is organized as follows. In Section \ref{sec:thoery} we present some needed notation, definitions and simple results from interpolation theory and algebra. Next, in Section \ref{ROFModelOnTheGraph} the TV minimization problem on the graph is formulated, the proposed algorithm is given and its convergence is proved. Thereafter, in Section \ref{SecNumExpts}, we present numerical experiments in order to compare the performance of the proposed algorithm with other iterative TV denoising algorithms and the BM3D image denoising method. Finally, discussion and concluding remarks are given in Section \ref{SecConclznz}.
\section{Notation and definitions}\label{sec:thoery}
\par In this section we briefly introduce the necessary mathematical theory and notation needed for presentation of the proposed algorithm. 
\subsection{Interpolation theory}\label{Interpolationtheori}
Let $X_0$ and $X_1$ be two Banach spaces. They form a Banach couple $\left(X_0,X_1\right)$ if there exists a Hausdorff topological vector space $\mathcal{H}$ in which both $X_0$ and $X_1$ are linearly and continuously embedded. For an introduction to the theory of interpolation, we refer to the book \cite{JBerghJL}. When $\left(X_0,X_1\right)$ is a Banach couple, then the sum $X_{0}+X_{1}$ given by $X_{0}+X_{1}=\{x\in\mathcal{H}:\ x=x_{0}+x_{1},\ x_{j}\in X_{j},\ j=0,1\}$ is well defined, and can be shown to be a Banach space under the norm $\left\|x\right\|_{X_{0}+X_{1}}=\inf\left(\left\|x_{0}\right\|_{X_{0}}+\left\|x_{1}\right\|_{X_{1}},\ x=x_{0}+x_{1},\ x_{j}\in X_{j},\ j=0,1\right)$. Furthermore, given a Banach couple $\left(X_0,X_1\right)$, an element $u_{0}\in X_0+X_1$ and a positive parameter $t$ the Peetre's $K$-functional is defined by 
\begin{equation*}
	K\left(t,u_{0};X_0,X_1\right) = \inf_{u\in X_1}\left(\left\|u_{0}-u\right\|_{X_0}+t\left\|u\right\|_{X_1}\right).
\end{equation*}
The $K$-functional is very important for the so-called $K$-method of real interpolation which generates families of real interpolation spaces between $X_{0}$ and $X_{1}$. The $K$-functional is a particular case of the more general $L$-functional which, for given $1\leq p_0,p_1<\infty$, is defined by
\begin{equation} \label{Lfunctional}
	L_{p_0,p_1}\left(t,u_{0};X_0,X_1\right) = \inf_{u\in X_1}\left(\frac{1}{p_{0}}\left\|u_{0}-u\right\|^{p_0}_{X_0}+\frac{t}{p_{1}}\left\|u\right\|^{p_1}_{X_1}\right).
\end{equation}
We need the following definitions of exact minimizers and optimal decomposition.
\vskip2mm 
\begin{definition}[\textbf{Exact minimizers}]
\label{NearAndExactMinimizers}
\par We say that the element $u_{opt}\in X_1$ is an \textit{exact minimizer} for the functional \eqref{Lfunctional} if
\begin{align*}
	\frac{1}{p_{0}}\left\|u_{0}-u_{opt}\right\|^{p_0}_{X_0}+\frac{t}{p_{1}}\left\|u_{opt}\right\|^{p_1}_{X_1}=L_{p_0,p_1}\left(t,u_{0};X_0,X_1\right).
\end{align*}
\end{definition}
\begin{definition}[\textbf{Optimal decomposition}]
\par If $u_{opt}\in X_1$ is an \textit{exact minimizer} for \eqref{Lfunctional}, then we call $u_{0} = u_{opt} + \left(u_{0}-u_{opt}\right)$ an \textit{optimal decomposition} for \eqref{Lfunctional}.
\end{definition}
\vskip2mm
\begin{remark}
\par It is important to note that an exact minimizer, and therefore an optimal decomposition, does not always exist.
\end{remark}
\vskip2mm
The $L_{p_{0},p_{1}}$-functional appears in regularization of inverse problems where the second term in the expression \eqref{Lfunctional} is called a penalty term or regularization term. Note that the total variation regularization functional \eqref{TheMinProblem} above is a particular case of the $L$-functional \eqref{Lfunctional} for $p_0=2$, $p_1=1$ and for the spaces $X_0=L^{2}\left(\Omega\right)$ and $X_1=BV\left(\Omega\right)$.
\subsection{Some algebra}\label{sec:algebra}
We start with the definition of the notion of annihilator.
\begin{definition}
\par Let $X$ be a Banach space and let $Z$ be a subspace of $X$. The \textit{annihilator} of $Z$ denoted $\Ann\left(Z\right)$ is the set of bounded linear functionals that vanish on $Z$. That is the set defined by
\begin{align*}
	\Ann\left(Z\right) = \left\{ x^{\ast}\in X^{\ast}:\ \left\langle x^{\ast},z\right\rangle = 0,\ \textrm{for all}\ z\in Z \right\},
\end{align*}
where $X^{\ast}$ is the dual space of $X$ and $\left\langle x^{\ast},z\right\rangle$ denotes the action of the bounded linear functional $x^{\ast}\in X^{\ast}$ on the element $z\in Z$.
\end{definition}
\vskip2mm
We will make use of the following result in the sequel. 
\vskip2mm
\begin{lemma}
\label{TheLemma}
\par Let $X$ be a Banach space with dual space $X^{\ast}$, $x_{0}\in X$ and let $Z$ be a finite-dimensional subspace of $X$. Then
\begin{align*}
	\inf_{z\in Z}\left\|x_{0}+z\right\|_{X} = \sup_{x^{\ast}\in\mathcal{B}_{X^{\ast}}\cap \Ann\left(Z\right)}\left\langle x^{\ast},x_{0}\right\rangle,
\end{align*}
where $\mathcal{B}_{X^{\ast}}$ is the unit ball of $X^{\ast}$.
\end{lemma}
\vskip2mm
\begin{proof}
The case $x_{0}\in Z$ is obvious. From now on we suppose $x_{0}\notin Z$. Let us take an arbitrary $x^{\ast}\in \mathcal{B}_{X^{\ast}}\cap \Ann\left(Z\right)$. Then we have
\begin{align*}
  \left\langle x^{\ast},x_{0}\right\rangle = \left\langle x^{\ast},x_{0}+z\right\rangle\leq \left\|x^{\ast}\right\|_{X^{\ast}}\left\|x_0+z\right\|_{X}\leq \left\|x_0+z\right\|_{X},\ \forall z\in Z.
\end{align*}
Therefore since $x^{\ast}\in\mathcal{B}_{X^{\ast}}\cap \Ann\left(Z\right)$ and $z\in Z$ are arbitrary, we have that
\begin{align}
\label{AnnihofTVVDistLOn1}
\inf_{z\in Z}\left\|x_{0}+z\right\|_{X} \geq \sup_{x^{\ast}\in\mathcal{B}_{X^{\ast}}\cap \Ann\left(Z\right)}\left\langle x^{\ast},x_{0}\right\rangle.
\end{align}
In order to prove the reverse inequality, let us consider the space $W$, which is the algebraic sum of the span of $x_0$ and the space $Z$:
\begin{align*}
	W = \span \left\{x_{0}\right\} + Z = \left\{w\in X:\ w=\lambda x_{0}+z,\ z\in Z \ \textrm{and}\ \lambda \in \mathbb{R}\right\},
\end{align*}
and take $z_{0}\in Z$ such that $\inf_{z\in Z}\left\|x_{0}+z\right\|_{X} = \left\|x_{0}+z_{0}\right\|_{X}$. The existence of such $z_{0}$ follows from the assumption that $Z$ is a finite-dimensional subspace of $X$. Without loss of generality we can assume that $\left\|x_{0}+z_{0}\right\|_{X} = 1$. Since $W$ is a normed vector space, it is possible to consider its dual space. Further, as $Z$ is a linear subspace of $W$, $x_{0}+z_{0}\in W$ and $x_{0}+z_{0}\notin Z$, the Hahn-Banach Theorem (see for example Corollary II.3.13 in \cite{Dunford1}) gives that there exists a bounded linear functional $x^{\ast}_{0}\in W^{\ast}$ such that
\begin{align*}
\left\langle x^{\ast}_{0},z\right\rangle = 0\ \textrm{for all}\ z\in Z\ \textrm{and}\ \left\langle x^{\ast}_{0},x_{0}+z_0\right\rangle = 1. 
\end{align*}
It follows that 
\begin{align}
\label{ThsyWEhavepdigt}
 x^{\ast}_{0} \in \Ann\left(Z\right)\ \textrm{and}\ \left\langle x^{\ast}_{0},x_{0}\right\rangle = 1.
\end{align}
Let us now investigate the action of $x^{\ast}_{0}$ on $W$. Let $w=\lambda x_0+z$ be an element of $W$ for some $\lambda\in\mathbb{R}$ and $z\in Z$. Then we have
\begin{align} \label{AnnihofTVVDistLOn5}
  \left\langle x^{\ast}_{0},w\right\rangle&=\left\langle x^{\ast}_{0},\lambda x_{0}+z\right\rangle = \left\langle x^{\ast}_{0},\lambda x_{0}+ \lambda z_{0} +z - \lambda z_{0}\right\rangle \nonumber \\
&= \lambda \left\langle x^{\ast}_{0},x_{0}+ z_{0} \right\rangle + \left\langle x^{\ast}_{0},z - \lambda z_{0}\right\rangle=\lambda,
\end{align}
because $\left\langle x^{\ast}_{0},x_{0}+z_0\right\rangle = 1$ and $\left\langle x^{\ast}_{0},z - \lambda z_{0}\right\rangle=0$ since $x^{\ast}_{0} \in \Ann\left(Z\right)$ and $ z-\lambda z_{0}\in Z$. Let us now describe the unit ball $\mathcal{B}_{W}$ of $W$. Suppose that $w = \lambda x_{0}+z \in \mathcal{B}_{W}$ where $\lambda\neq 0$. We have that
\begin{align}
\label{AnnihofTVVDistLOn6}
	1\geq\left\|w\right\|_{X}=\left\|\lambda x_{0}+z\right\|_{X} = \left|\lambda\right|\left\|x_{0}+\frac{z}{\lambda}\right\|_{X}\geq \left|\lambda\right|\left\|x_{0}+z_{0}\right\|_{X} = \left|\lambda\right|.
\end{align}
Therefore $w = \lambda x_{0}+z \in \mathcal{B}_{W}$ implies that $\left|\lambda\right|\leq 1$. 
From \eqref{AnnihofTVVDistLOn5} and \eqref{AnnihofTVVDistLOn6}, it follows that 
\begin{align*}
	\left\|x^{\ast}_{0}\right\|_{W^{\ast}}= \sup_{w\in\mathcal{B}_{W}}\left\langle x^{\ast}_{0},w\right\rangle = \sup_{w\in\mathcal{B}_{W}} \lambda = \sup_{\left|\lambda\right|\leq 1} \lambda =1.
\end{align*}
By invoking the Hahn-Banach theorem (see for example Theorem II.3.11 in \cite{Dunford1}), we can extend the functional $x^{\ast}_{0}$ to a functional $\widetilde{x^{\ast}_{0}}\in X^{\ast}$ such that $\widetilde{x^{\ast}_{0}}\rvert_{W}=x^{\ast}_{0}$ and $\left\|\widetilde{x^{\ast}_{0}}\right\|_{X^{\ast}} = \left\|x^{\ast}_{0}\right\|_{W^{\ast}} = 1$.
From this and \eqref{ThsyWEhavepdigt} we conclude that $\widetilde{x^{\ast}_{0}}\in \mathcal{B}_{X^{\ast}}\cap \Ann\left(Z\right)$.
 It follows that
\begin{align} \label{AnnihofTVVDistg2}
	\inf_{z\in Z}\left\|x_{0}+z\right\|_{X}&=\left\|x_{0}+z_{0}\right\|_{X}=\langle x^{\ast}_{0},x_{0}+z_{0}\rangle=\langle\widetilde{x^{\ast}_{0}},x_0+z_0\rangle \nonumber\\
	&\leq \sup_{x^{\ast}\in \mathcal{B}_{X^{\ast}}\cap \Ann\left(Z\right)}\langle x^{\ast},x_0+z_0\rangle=\sup_{x^{\ast}\in \mathcal{B}_{X^{\ast}}\cap \Ann\left(Z\right)}\langle x^{\ast},x_0\rangle.
\end{align}
Putting \eqref{AnnihofTVVDistLOn1} and \eqref{AnnihofTVVDistg2} together, we obtain
\begin{align*}
	\inf_{z\in Z}\left\|x_{0}+z\right\|_{X} = \sup_{x^{\ast}\in \mathcal{B}_{X^{\ast}}\cap \Ann\left(Z\right)}\left\langle x^{\ast},x_0 \right\rangle
\end{align*}
which concludes the proof.
\end{proof}
\vskip2mm
\section{Introducing iterative TV minimization on the graph}\label{ROFModelOnTheGraph}
\subsection{A graph specific problem formulation}
Suppose we have an observed noisy image $u_{0}\in L^{2}$ defined on the domain $\Omega=\left(0,1\right)^{2}\subset\mathbb{R}^{2}$ which is a degraded version of the original true image $u\in BV\left(\Omega\right)$ according to the linear model \eqref{eq:linmodel}. The ROF model suggests to take as an approximation to the original image $u$ the function $u_{opt}\in BV$ which is the exact minimizer for the $L_{2,1}$-functional of the couple $\left(L^{2},BV\right)$:
\begin{equation}
\label{ROFModelIntro}
L_{2,1}\left(  t,u_{0};L^{2},BV\right)  =\inf_{u\in BV}\left(  \frac{1}%
{2}\left\Vert u_{0}-u\right\Vert _{L^{2}}^{2}+t\left\Vert u\right\Vert
_{BV}\right)  ,\ \text{for some}\ t>0.
\end{equation}
We will use the following anisotropic $BV$ seminorm:
\begin{align*} 
\left\Vert u\right\Vert _{BV\left(  \Omega\right)  }=\int_{0}^{1}\operatorname{var}_{x}u\left(y\right)  dy+\int_{0}^{1}\operatorname{var}_{y}u\left(x\right)  dx,
\end{align*}
where
\[
\operatorname{var}_{x}u\left(y\right)  =\sup_{0\leq x_{1}\leq\ldots\leq x_{n}\leq1}%
\sum_{j=1}^{n-1}\left\vert u\left(  x_{j+1},y\right)  -u\left(  x_{j}%
,y\right)  \right\vert
\]
is the total variation of $u$ along the the horizontal axis for a given $y$, and
\[
\operatorname{var}_{y}u\left(x\right)  =\sup_{0\leq y_{1}\leq\ldots\leq y_{n}\leq1}%
\sum_{i=1}^{n-1}\left\vert u\left(  x,y_{i+1}\right)  -u\left(  x,y_{i}%
\right)  \right\vert
\]
is the total variation of $u$ along the vertical axis for a given $x$. The reason for choosing this $BV$ seminorm is that it suggests a convenient formulation of total variation in the graph setting, see \eqref{ROFMscGVEGDef} and \eqref{ROFMscGVEGDef2}.
\par We use a standard approach when discretizing the functional \eqref{ROFModelIntro}, i.e., we divide $\Omega$ into $N\times N$ square cells and instead of the space $L^{2}(\Omega)$ consider its finite-dimensional subspace $S_{N}$ consisting of functions that are piecewise constant on each cell. Throughout, we consider our discretization grid as a 2D Cartesian coordinate in \textit{screen space}, i.e., the same way matrices are represented on the computer. We define
\begin{align*}
	S_{N} = \left\{u= \sum^{N}_{i,j=1}u_{ij}\scalebox{1.5}{$\chi$}_{ij},\ \scalebox{1.5}{$\chi$}_{ij}\left(x,y\right) = \left\{
\begin{array}{rl}
1 & \text{if } \frac{j-1}{N}< x<\frac{j}{N}\ \textrm{and}\ \frac{i-1}{N}< y<\frac{i}{N}\\
\vspace{-4mm}\\
0 & \text{otherwise}.
\end{array} \right. \right\}.
\end{align*}
It is clear that the $BV$ seminorm of a function $u\in S_{N}$ is equal to%
\[
\left\Vert u\right\Vert _{BV\left(S_{N}\right)}=\frac{1}{N}\left(  \sum_{i=1}^{N}\sum
_{j=1}^{N-1}\left\vert u_{i,j+1}-u_{ij}\right\vert +\sum_{j=1}^{N}\sum
_{i=1}^{N-1}\left\vert u_{i+1,j}-u_{ij}\right\vert \right).
\]
Therefore the discrete analogue of the functional \eqref{ROFModelIntro} can be written as%
\begin{align}
\label{ROFModelIntroDisc}
L_{2,1}\left(t,u_{0};L^{2},BV\right)&=\inf_{u\in S_{N}}\Bigg(\frac
{1}{2N^{2}}\left(\sum_{i,j=1}^{N}\left(  u_{0_{ij}}-u_{ij}\right)
^{2}\right) \nonumber\\&+\frac{t}{N}\left(\sum_{i=1}^{N}\sum_{j=1}^{N-1}\left\vert
u_{i,j+1}-u_{ij}\right\vert +\sum_{j=1}^{N}\sum_{i=1}^{N-1}\left\vert
u_{i+1,j}-u_{ij}\right\vert \right)  \Bigg).
\end{align}
\subsubsection{Graph notations}
We now turn to the framework of graph which generalizes the problem \eqref{ROFModelIntroDisc}. Let $G=\left(V,E\right)$ be a finite, directed and connected\footnote{The analysis can be extended to disconnected graphs as the components of the graph are considered separately. For convenience, we have therefore chosen to only consider connected graphs in this paper.} graph with $N$ vertices $V=\left\{v_{1},v_{2},\ldots,v_{N}\right\}$ and $M$ directed edges
$E=\left\{e_{1},e_{2},\ldots,e_{M}\right\}$ where each edge is determined by a pair of vertices, i.e. $e_{k}=\left(v_{i},v_{j}\right)$ for some $i,j\in\left\{1,2,\ldots,N\right\}$ and $k=1,2,\ldots,M$. We assume that the edge $e_{k}=\left(v_{i},v_{j}\right)$ is directed from the vertex $v_{i}$ to the vertex $v_{j}$. Let $S_{V}=\left\{f:\ f:V\rightarrow\mathbb{R}\right\} $ denote the $N-$dimensional space of real-valued functions defined on $V$ and let $S_{E}=\left\{g:\ g:E\rightarrow \mathbb{R}\right\}$ denote the $M$-dimensional space of real-valued functions defined on $E$.
\par The \emph{gradient} operator $\operatorname{grad}:S_{V}\rightarrow S_{E}$ is defined by 
\begin{align*}
	\operatorname{grad}f(e)=f(v_{j})-f(v_{i}),e=(v_{i},v_{j})\in E.
\end{align*}
We define inner products on $S_{E}$ and $S_{V}$ according to
\begin{align*}
 \left\langle f_{1},f_{2}\right\rangle_{S_{E}}=\sum_{e\in E}f_{1}(e)f_{2}(e)
\end{align*}
and
\begin{align*}
 \left\langle g_{1},g_{2}\right\rangle_{S_{V}}=\sum_{v\in V}g_{1}(v)g_{2}(v).
\end{align*}
It is easy to show that the \emph{divergence} operator $\operatorname{div}:S_{E}\rightarrow S_{V}$ given by
\begin{align*}
\operatorname{div}g(v_{j})=\sum_{i:(v_{i},v_{j})\in E}g\left((v_{i},v_{j})\right)-\sum_{k:(v_{j},v_{k})\in
E}g\left((v_{j},v_{k})\right).
\end{align*}
is conjugate to $\operatorname{grad}$, i.e.
\begin{align}
\label{TheOtherLemma}
	\left\langle \operatorname{div}g,f\right\rangle _{S_{V}}=\left\langle
g,\operatorname{grad}f\right\rangle _{S_{E}},\ \forall f\in S_{V}, \forall g\in S_{E}. 
\end{align}
If we consider elements of $S_{E}$ as flows on the graph $G=\left(V,E\right)$, the divergence at a vertex can be interpreted as the difference between the total incoming flows and the total outgoing flows.
\par
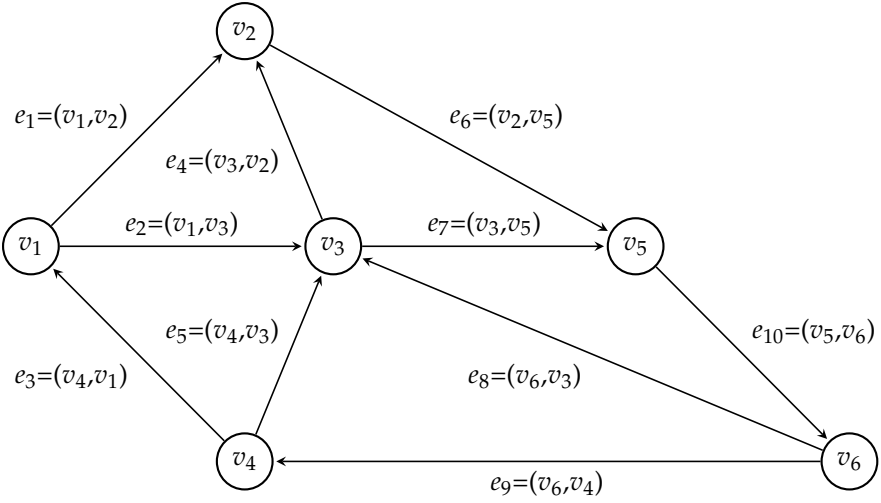
\begin{figure}[h]
  \begin{tikzpicture}[
    > = stealth, 
    shorten > = 1pt, 
    auto,
    node distance = 4cm, 
    semithick 
    ]
    \tikzstyle{every state}=[
    draw = black,
    thick,
    fill = white,
    minimum size = 4mm
    ]
    \node[state] (v1) {$v_1$};
    \node[state] (v2) [above right of=v1] {$v_2$};
    \node[state] (v3) [right of=v1] {$v_3$};
    \node[state] (v4) [below right of=v1] {$v_4$};
    \node[state] (v5) [right of=v3] {$v_5$};
    \node[state] (v6) [below right of=v5] {$v_6$};
    \path[->] (v1) edge node {$e_1$=($v_1$,$v_2$)} (v2);
    \path[->] (v1) edge node {$e_2$=($v_1$,$v_3$)} (v3);
    \path[->] (v4) edge node {$e_3$=($v_4$,$v_1$)} (v1);
    \path[->] (v3) edge node {$e_4$=($v_3$,$v_2$)} (v2);
    \path[->] (v4) edge node {$e_5$=($v_4$,$v_3$)} (v3);
    \path[->] (v2) edge node {$e_6$=($v_2$,$v_5$)} (v5);
    \path[->] (v3) edge node {$e_7$=($v_3$,$v_5$)} (v5);
    \path[->] (v6) edge node {$e_8$=($v_6$,$v_3$)} (v3);
    \path[->] (v6) edge node {$e_9$=($v_6$,$v_4$)} (v4);
    \path[->] (v5) edge node {$e_{10}$=($v_5$,$v_6$)} (v6);
  \end{tikzpicture}
  \caption{Graph illustrating the notation and node relations for the discrete gradient and divergence operators.}
  \label{fig:graph}
\end{figure}

 The graph in Figure \ref{fig:graph} illustrates the definition for the gradient and the divergence operator by an example for the case $N=6$ and $M=10$. For example, $\operatorname{grad}f(e_6)=f(v_{5})-f(v_{2})$ is the gradient at $e_6=(v_2$,$v_5$) and $\operatorname{div}g(v_{3})=g\left((v_{1},v_{3})\right)+g\left((v_{4},v_{3})\right)+g\left((v_{6},v_{3})\right)-g\left((v_{3},v_{2})\right)-g\left((v_{3},v_{5})\right)$ is the divergence at $v_{3}$.
\vskip2mm
\begin{remark}
\label{RemarkABTBVstar}
\par The operator $\operatorname{grad}$ has a kernel given by
\begin{align*}
	\ker\left(\operatorname{grad}\right) = \left\{ f\in S_{V}:\ f = C,\ \textrm{for some}\ C\in\mathbb{R} \right\},
\end{align*}
and its orthogonal complement coincides with its annihilator and is given by
\begin{align}
\label{AAnnihAGen010}
	\left(\ker\left(\operatorname{grad}\right)\right)^{\bot} = \Ann\left(\ker\left(\operatorname{grad}\right)\right) = \left\{ F\in S_{V}:\ \sum_{v\in V}F\left(v\right)=0 \right\}.
\end{align}
Since $\operatorname{div}$ is the conjugate operator of $\operatorname{grad}$, the fundamental theorem of linear algebra ensures that
\begin{align}
\label{AKernelAstarGen}
	\operatorname{im} \left(\operatorname{div}\right) = \Ann \left( \ker \left(\operatorname{grad}\right)  \right)\ \textrm{and}\ \operatorname{im} \left(\operatorname{grad}\right) = \Ann \left( \ker \left(\operatorname{div}\right)  \right)
\end{align}
where $\operatorname{im}\left(A\right)$ denotes the image of the operator $A$.
\end{remark}
\vskip2mm
\par An observed image $u_{0}\in S_{N}$ can be considered as an element of $S_{V}$ for a graph $G=(V,E)$ where the cells are represented by the vertices in $V$ and pairs of adjacent cells are represented by edges in $E$ (any direction of the edges can be chosen). The functional \eqref{ROFModelIntroDisc} can then be written as%
\begin{equation}
\label{ROFModelIntroDiscGVE}
L_{2,1}\left(  t,u_{0};L^{2},BV\right)  =\inf_{u\in S_{V}}\left(\frac
{1}{2N^{2}}\left\Vert u_{0}-u\right\Vert _{\ell^{2}(S_{V})}^{2}+\frac{t}%
{N}\left\Vert \operatorname{grad}u\right\Vert _{\ell^{1}(S_{E})}\right).
\end{equation}
It is clear that the exact minimizer of \eqref{ROFModelIntroDiscGVE} coincides with the exact minimizer of%
\begin{align*}
L_{2,1}\left(s,u_{0};\ell^{2}(S_{V}),BV(S_{V})\right)=\inf_{u\in S_{V}%
}\left(  \frac{1}{2}\left\Vert u_{0}-u\right\Vert _{\ell^{2}(S_{V})}%
^{2}+s\left\Vert \operatorname{grad}u\right\Vert _{\ell^{1}(S_{E})}\right),\ s=Nt.
\end{align*}
This observation leads to the following analogue of the ROF model on a general finite, connected and directed
graph.
\vskip2mm
\begin{problem}
\label{TheProblemRG}
\par Suppose that we know the function $u_{0}\in S_{V}$. For given $t>0$, find the exact
minimizer of the functional%
\begin{equation*}
L_{2,1}\left(  t,u_{0};\ell^{2}(S_{V}),BV(S_{V})\right)  =\inf_{u\in
BV(S_{V})}\left(  \frac{1}{2}\left\Vert u_{0}-u\right\Vert _{\ell^{2}(S_{V}%
)}^{2}+t\left\Vert u\right\Vert _{BV(S_{V})}\right),
\end{equation*}
where
\begin{subequations}
\begin{align}\label{ROFMscGVEGDef}
\left\Vert u\right\Vert _{\ell^{2}(S_{V})}&=\left(  \sum_{v\in V}\left(
u(v)\right)  ^{2}\right)  ^{\frac{1}{2}},
\quad
\left\Vert u\right\Vert _{BV(S_{V}%
)}=\left\Vert \operatorname{grad}u\right\Vert _{\ell^{1}(S_{E})}, \\
\label{ROFMscGVEGDef2}
& \mbox{ and } \qquad  \left\Vert
\psi\right\Vert _{\ell^{1}(S_{E})}=\sum_{e\in E}\left\vert \psi(e)\right\vert.
\end{align}
\end{subequations}
\end{problem}
\subsubsection{Description of the ball of dual space to \texorpdfstring{$BV(S_{V})$}{Lg}}
In order to describe our algorithm for Problem \ref{TheProblemRG}, we first need a description of the ball of the dual space to $BV(S_{V})$.
\par
\begin{figure}[t]
\centering
\hspace{-20mm}
\begin{overpic}[width=0.8\textwidth]{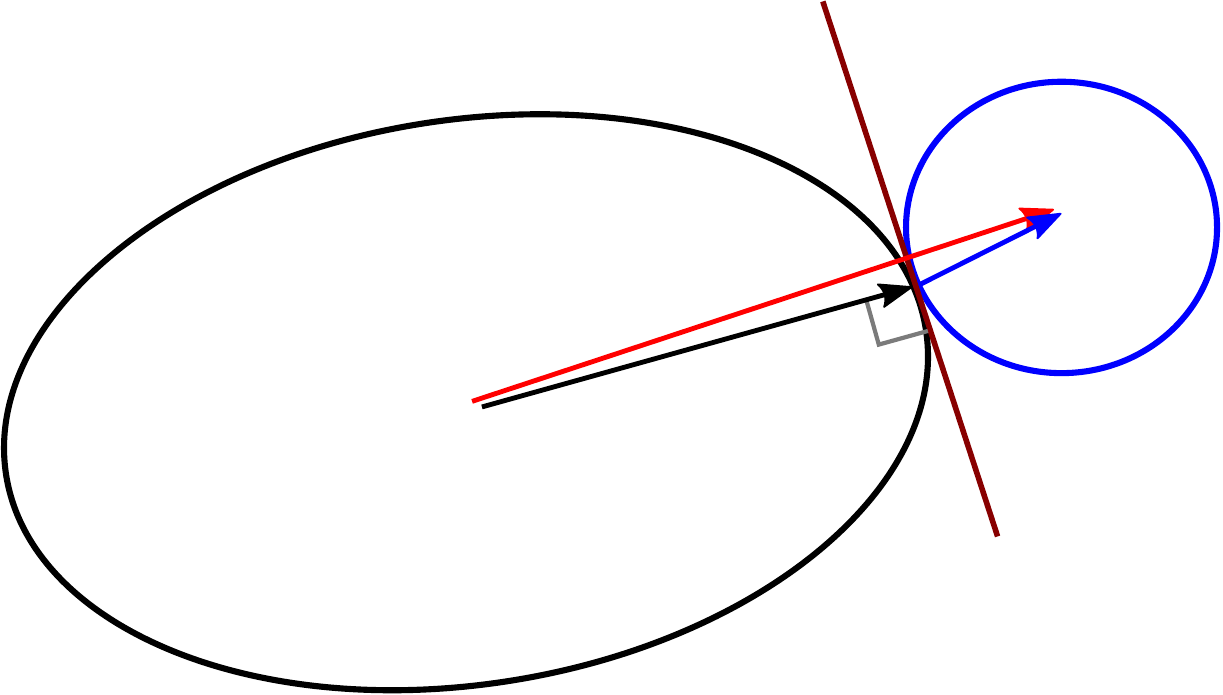} 
 \put (63,54) {$\displaystyle T$}
 \put (36,22) {$\displaystyle O$}
 \put (55,25) {$\displaystyle u_0 - u_{opt}$}
 \put (80,32) {$\displaystyle u_{opt}$}
 \put (88,37) {$\displaystyle u_{0}$}
 \put (10,50) {$\displaystyle t \mathcal{B}_{X^*} = \{ \psi \in \mathbb{R}^n : || \psi ||_{X^*} \leq t \}$}
 \put (70,52) {$\displaystyle  \{ u \in \mathbb{R}^n : || u_0 - u ||_{\ell^2} \leq || u_{opt} ||_{\ell^2} \}$}
\end{overpic}
\caption{Geometrical illustration of the ball of dual space
and the position of the element of best approximation $u_{opt}=u_{0}-\underset{\psi\in 
t\mathcal{B}_{X^{\ast}}}\argmin \left\|u_{0}-\psi\right\|_{\ell^{2}}$. The hyperplane $T$ is orthogonal to $u_{0}-u_{opt}$.}
\label{fig:GeoIllutrAppr}
\end{figure}
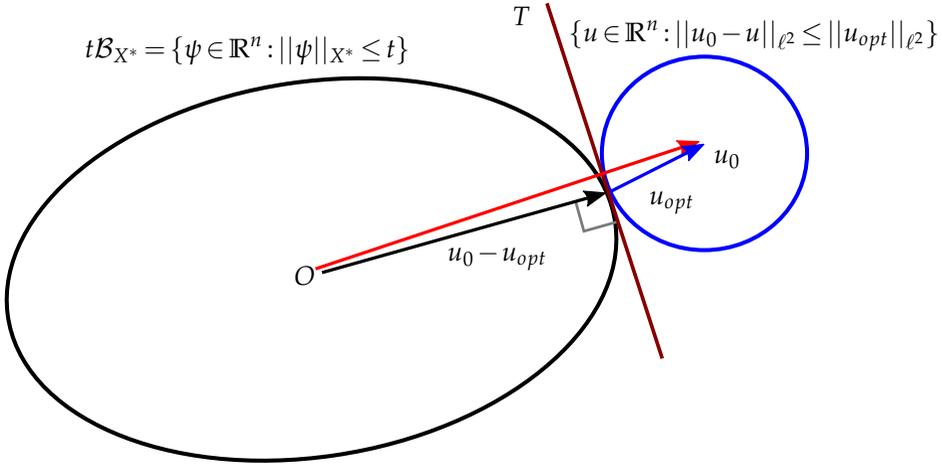

 It was shown in \cite{NatanJaphet2} that the exact minimizer $u_{opt}$ for the $L_{2,1}$- functional for the couple $\left(\ell^{2},X\right)$, where $X$ is a Banach space,
\[
L_{2,1}\left(  t,u_{0};\ell^{2},X\right)  =\inf_{u\in X}\left(  \frac{1}%
{2}\left\Vert u_{0}-u\right\Vert _{\ell^{2}}^{2}+t\left\Vert u\right\Vert
_{X}\right),
\]
is equal to the difference between $u_{0}$ and the nearest element to $u_{0}$
of the ball of radius $t>0$ of the space $X^{\ast}$, i.e., $u_{opt}=u_{0}-\underset{\psi\in t\mathcal{B}_{X^{\ast}}}\argmin \left\|u_{0}-\psi\right\|_{\ell^{2}}$. Figure \ref{fig:GeoIllutrAppr} provides a geometrical illustration of the optimal decomposition.
\par Consider now $X=BV$. As $\left\|\cdot\right\|_{BV(S_{V})}$ is a seminorm on $S_{V }$, we restrict to the subspace $\left(\ker\left(\operatorname{grad}\right)\right)^{\bot}$ where $\left\|\cdot\right\|_{BV(S_{V})}$ is a norm.
The dual space $BV^{\ast}(S_{V})$ is then $\left(\ker\left(\operatorname{grad}\right)\right)^{\bot}$ equipped with the norm defined by%
\begin{align}
\label{BVstarBallDefnm}
	\left\Vert \psi\right\Vert _{BV^{\ast}(S_{V})}=\sup_{\left\Vert h\right\Vert
_{BV(S_{V})}\leq 1}\left\langle \psi,h\right\rangle_{S_{V}}.
\end{align}
We have the following characterization of the unit ball of $BV^{\ast}(S_{V})$:
\vskip2mm
\begin{theorem}
\label{BVstarBall}
\par The unit ball of the space $BV^{\ast}(S_{V})$ is equal to the image of the unit 
ball of the space $\ell^{\infty}(S_{E})$ under the operator $\operatorname{div}$, i.e.,%
\begin{align*}
	\mathcal{B}_{BV^{\ast}(S_{V})}=\operatorname{div}\left(\mathcal{B}_{\ell^{\infty}(S_{E})}\right).
\end{align*}
\end{theorem}
\vskip2mm
\begin{proof}
Let us consider an arbitrary $\psi\in BV^{\ast}\left(S_{V}\right)$. From relations \eqref{AAnnihAGen010} and \eqref{AKernelAstarGen}, we conclude that $BV^{\ast}\left(S_{V}\right) = \operatorname{im} \left(\operatorname{div}\right)$. Therefore for all $\psi\in BV^{\ast}\left(S_{V}\right)$, there exists at least one $g\in S_{E}$ such that $\psi = \operatorname{div}g$. Fix $g_{0}\in S_{E}$ such that $\psi = \operatorname{div}g_{0}$. We have
\begin{align}
\label{ElemLinAlgebra}
	\underset{\psi=\operatorname{div}g}\inf\left\|g\right\|_{\ell^{\infty}\left(S_{E}\right)} = \underset{\varphi\in\ker \left(\operatorname{div}\right)}\inf\left\|g_{0}+\varphi\right\|_{\ell^{\infty}\left(S_{E}\right)}.
\end{align}
By applying Lemma \ref{TheLemma} and using expression \eqref{AKernelAstarGen}, together with \eqref{TheOtherLemma}, $\psi=\operatorname{div}g_{0}$, \eqref{ROFMscGVEGDef} and \eqref{BVstarBallDefnm} we derive
\begin{align*}
\inf_{\psi=\operatorname{div}g}\left\|g\right\|_{\ell^{\infty}\left(S_{E}\right)} &= \underset{\varphi\in\ker \left(\operatorname{div}\right)}\inf\left\|g_{0}+\varphi\right\|_{\ell^{\infty}\left(S_{E}\right)}=\underset{f\in\mathcal{B}_{\ell^{1}\left(S_{E}\right)}\cap \operatorname{im} \left(\operatorname{grad}\right)}\sup\left\langle f,g_{0}\right\rangle_{S_{E}}
\\
& = \sup_{\left\|\operatorname{grad}h\right\|_{\ell^{1}\left(S_{E}\right)}\leq 1} \left\langle \operatorname{grad}h,g_{0}\right\rangle_{S_{E}}=\sup_{\left\|\operatorname{grad}h\right\|_{\ell^{1}\left(S_{E}\right)}\leq 1} \left\langle h,\operatorname{div}g_{0}\right\rangle_{S_{V}}\\
&= \sup_{\left\|h\right\|_{BV\left(S_{V}\right)}\leq 1} \left\langle h,\psi\right\rangle_{S_{V}}= \left\|\psi\right\|_{BV^{\ast}\left(S_{V}\right)}.
\end{align*}
From this follows that
\begin{align*}
	\left\|\psi\right\|_{BV^{\ast}\left(S_{V}\right)}\leq 1\ \textrm{if and only if}\ \inf_{\psi=\operatorname{div}g}\left\|g\right\|_{\ell^{\infty}\left(S_{E}\right)}\leq 1.
\end{align*}
So, it is clear that $BV^{\ast}\left(S_{V}\right)\supset\operatorname{div}\left(\mathcal{B}_{\ell^{\infty}\left(S_{E}\right)} \right)$. Note next that the infimum in \eqref{ElemLinAlgebra} is attained because $\ker \left(\operatorname{div}\right)$ is a subspace of the finite-dimensional space $S_{E}$. Therefore, for each $\psi\in\mathcal{B}_{BV^{\ast}\left(S_{V}\right)}$ there exists an element
 \begin{align*}
	 g_{\psi} \in g_{0}+ \ker \left(\operatorname{div}\right),\ \textrm{such that}\ \left\|g_{\psi}\right\|_{\ell^{\infty}\left(S_{E}\right)}\leq 1\ \textrm{and}\ \operatorname{div}g_{\psi}=\psi.
 \end{align*}
We conclude that 
\begin{align*}
	\mathcal{B}_{BV^{\ast}\left(S_{V}\right)} = \operatorname{div}\left(\mathcal{B}_{\ell^{\infty}\left(S_{E}\right)} \right).
\end{align*}
\end{proof}
\vskip2mm
\subsection{Algorithm}\label{sec:Algorithm}
Algorithm \ref{ALGO} below embodies our algorithmic contribution for computing the ROF-minimizer $u_{opt}$ and we will now describe its construction in detail.
\par The core of the algorithm is the construction of the element $\widetilde{\psi} = \left(u_{0}-u_{opt}\right)\in t\mathcal{B}_{BV^{\ast}(S_{V})}$ that satisfies
\begin{align*}
	\left\|u_{0}-\widetilde{\psi}\right\|_{\ell^{2}\left(S_{V}\right)}=\inf_{\psi\in  t\mathcal{B}_{BV^{\ast}(S_{V})}} \left\|u_{0}-\psi\right\|_{\ell^{2}\left(S_{V}\right)}.
\end{align*}
From Theorem \ref{BVstarBall}, this is equivalent to construct a flow $g_{\widetilde{\psi}}\in t \mathcal{B}_{\ell^{\infty}\left(S_{E}\right)}$ such that 
\begin{align*}
	\left\|u_{0}-\operatorname{div}g_{\widetilde{\psi}}\right\|_{\ell^{2}\left(S_{V}\right)}=\inf_{g\in t \mathcal{B}_{\ell^{\infty}\left(S_{E}\right)}} \left\|u_{0}-\operatorname{div}g\right\|_{\ell^{2}\left(S_{V}\right)},
\end{align*}
and put $\widetilde{\psi}=\operatorname{div}g_{\widetilde{\psi}}$. Once this is done, $u_{opt} = u_{0}-\widetilde{\psi}$.
\par We now describe the steps of the algorithm in detail. Let $u_{0}$ be defined on $G=(V,E)$ with vertex set $V=\{v_1,\ldots,v_N\}$ and edge set $E=\{e_1,\ldots,e_M\}$. The parameter $t$ denotes a regularization parameter and $N_{iter}$ denotes the maximum number of iterations. The edge set is specifically defined as 
\begin{displaymath}
  e_{k} = \left( v_{i}, v_{j} \right)\in E,\ k=1,2,\ldots,M;\ \textrm{for some}\ i,j\in \left\{ 1,2,\ldots,N \right\}.
\end{displaymath}
Introduce the operator $T:t \mathcal{B}_{\ell^{\infty}\left(S_{E}\right)} \rightarrow t \mathcal{B}_{\ell^{\infty}\left(S_{E}\right)}$ given by
\begin{align} \label{U0}
T = T_{M}T_{M-1}T_{M-2}\ldots T_{2}T_{1}
\end{align}
where for $k=1,2,\ldots,M$, the operator $T_{k}: t \mathcal{B}_{\ell^{\infty}\left(S_{E}\right)} \rightarrow t \mathcal{B}_{\ell^{\infty}\left(S_{E}\right)}$ is defined as follows 
\begin{align} \label{U1}
T_{k}g\left(e\right) = \left\{
\begin{array}{lr}
\left\{
\begin{array}{rl}
Kg\left(e_{k}\right) & \text{,if } Kg\left(e_{k}\right)\in \left[ -t,+t \right];\\
-t & \text{,if } Kg\left(e_{k}\right)< -t ;\\
+t & \text{,if } Kg\left(e_{k}\right)> +t.
\end{array} \right.& \text{, if } e=e_k ;
\vspace{.0in}\\
g\left(e\right) & \text{, if } e\neq e_k.
\end{array} \right.
\end{align}
Here
\begin{align*} 
Kg\left(e_{k}\right) = \frac{\left[u_{0}(v_{j})-\operatorname{div}\backslash_{e_k}g(v_{j})\right]-\left[u_{0}(v_{i})-\operatorname{div}\backslash_{e_k}g(v_{i})\right]}{2}
\end{align*}
and
\begin{align*} 
\left\{
\begin{array}{lr}
\operatorname{div}_{\backslash e_{k}}g\left(v_{i}\right) = \operatorname{div}g\left(v_{i}\right) + g\left(e_{k}\right);\\
\operatorname{div}_{\backslash e_{k}}g\left(v_{j}\right) = \operatorname{div}g\left(v_{j}\right) - g\left(e_{k}\right);\\
\operatorname{div}_{\backslash e_{k}}g\left(v_{\ell}\right) = \operatorname{div}g\left(v_{\ell}\right),\ \forall \ell\neq i,j,
\end{array} \right.
\end{align*}
i.e. $\operatorname{div}_{\backslash e_{k}}$ is the divergence operator $\operatorname{div}$ without taking into account the flow on the edge $e_{k}$.
\par The constructed operator $T$ depends on the enumeration of the edges in $E$. However the results concerning $T$, i.e. Proposition \ref{TheOpTproperties} and Theorem \ref{PropoTconv} below, hold regardless of the specific enumeration of the edges. We will now point out a certain construction of $T$ which leads to a version of Algorithm \ref{ALGO} suitable for parallel computer architectures.
\par Colour the set of edges $E$ such that incident edges, i.e. edges that share a common vertex, have different colours. Denote by $E_{1},...,E_{L}$ the resulting disjoint subsets of $E$, $\cup_{i\in\left\{1,...,L\right\}}E_{i}=E$, from such a colouring with usage of $L$ different colours. Let $e_{i,1},....,e_{i,M_{i}}$ denote the edges of $E_{i}$ and define $T_{E_{i}}=T_{i,M_{i}}T_{i,M_{i}-1}...T_{i,1}$ where
\begin{align*}
T_{i,k}g\left(e\right) = \left\{
\begin{array}{lr}
\left\{
\begin{array}{rl}
Kg\left(e_{i,k}\right) & \text{,if } Kg\left(e_{i,k}\right)\in \left[ -t,+t \right];\\
-t & \text{,if } Kg\left(e_{i,k}\right)< -t ;\\
+t & \text{,if } Kg\left(e_{i,k}\right)> +t.
\end{array} \right.& \text{, if } e=e_{i,k} ;
\vspace{.0in}\\
g\left(e\right) & \text{, if } e\neq e_{i,k}.
\end{array} \right.
\end{align*}
Because the edges of $E_{i}$ are non-incident, it follows that the applications of $T_{i,k}$, $k=1,...,M_{i}$, can be done in arbitrary order without affecting the resulting update $T_{E_{i}}g$ of $g$. The associated computations can therefore be done in parallel. With $T_{E_{i}}$, $i=1,...,L$, given, the operator $T$ is then constructed according to
\begin{equation*}
T=T_{E_{L}}T_{E_{L-1}}...T_{E_{1}}.
\end{equation*}
\begin{algorithm}[!h]
\caption{\textbf{: ROF model on the graph}}
\label{ALGO}
\begin{algorithmic}
\STATE $n \leftarrow 0$
\STATE choose initial $g\in t \mathcal{B}_{\ell^{\infty}\left(S_{E}\right)}$
\STATE $g_{n} \leftarrow g$
\WHILE{$n < N_{iter}$}
\IF{$Tg_{n} = g_{n}$}
\STATE \textbf{stop}
\ELSE[$Tg_{n} \neq g_{n}$]
\STATE
Compute $g_{n+1} = Tg_{n}$ 
\ENDIF
\STATE $n=n+1$
\ENDWHILE
\STATE
Compute $\widetilde{\psi}=\operatorname{div}(g_{n})$ 
\RETURN $u_{opt} = u_{0}-\widetilde{\psi}$
\end{algorithmic}
\end{algorithm}
\subsection{Convergence results}
For Algorithm \ref{ALGO}, convergence is established in Theorem \ref{PropoTconv}. Its proof is based on the following proposition.
\vskip2mm
\begin{proposition}
\label{TheOpTproperties}
\par The operator $T:t \mathcal{B}_{\ell^{\infty}\left(S_{E}\right)} \rightarrow t \mathcal{B}_{\ell^{\infty}\left(S_{E}\right)}$ given by \eqref{U0}-\eqref{U1} is continuous and satisfies the following two conditions
\begin{itemize}
	\item[(1)] For any $g\in t \mathcal{B}_{\ell^{\infty}\left(S_{E}\right)}$, $\operatorname{div}g=\widetilde{\psi}\ \textrm{if and only if}\ Tg=g$;
	\vskip1mm
	\item[(2)] For any $g\in t \mathcal{B}_{\ell^{\infty}\left(S_{E}\right)}$, $\textrm{if}\ \operatorname{div}g\neq\widetilde{\psi}$ then $\left\| u_{0}-\operatorname{div}\left(Tg\right) \right\|_{\ell^{2}\left(S_{V}\right)} < \left\| u_{0}-\operatorname{div}g \right\|_{\ell^{2}\left(S_{V}\right)}$. 
\end{itemize}
\end{proposition}
\vskip2mm
\begin{proof}
Each operator $T_{k}$ is continuous because by definition, it is clear that small changes of $g\in t\mathcal{B}_{\ell^{\infty }(S_{E})}$ leads to small changes of $T_{k}$ and therefore $T$ is continuous as a product of continuous operators.
\par We now prove condition (1). Let $g\in t \mathcal{B}_{\ell^{\infty}\left(S_{E}\right)}$ and assume that $\operatorname{div}g=\widetilde{\psi}$. Take $e_{k}=(v_{i},v_{j})\in E$. We note that $u(e_{k})$ appears only in
the following two terms of $\left\Vert u_{0}-\operatorname{div}g\right\Vert_{\ell^{2}(S_{V})}^{2}$:
\begin{align}
\label{MIniThsTwo}
\left[ u_{0}(v_{j})-\operatorname{div}g(v_{j})\right]^{2}+\left[ u_{0}(v_{i})-
\operatorname{div}g(v_{i})\right] ^{2} \nonumber\\
=\left[u_{0}(v_{j})-\left(\operatorname{div}_{\backslash
e_{k}}g(v_{j})+g(e_{k})\right)\right] ^{2} + 
\left[u_{0}(v_{i})-\left(\operatorname{div}%
_{\backslash e_{k}}g(v_{i})-g(e_{k})\right)\right] ^{2}.
\end{align}%
Since $\operatorname{div}g=\widetilde{\psi}$, $g(e_{k})$ in particular must minimize \eqref{MIniThsTwo} 
in the interval $[-t,t]$. By Jensen's inequality we note that%
\begin{align*}
\xi (g(e_{k})) &=\left[u_{0}(v_{j})-\left(\operatorname{div}_{\backslash
e_{k}}g(v_{j})+g(e_{k})\right)\right]^{2}+ \left[u_{0}(v_{i})-\left(\operatorname{div}%
_{\backslash e_{k}}g(v_{i})-g(e_{k})\right)\right] ^{2}\\
&\geq 2\left(\frac{\left[u_{0}(v_{j})-\operatorname{div}_{\backslash
e_{k}}g(v_{j})\right]+\left[u_{0}(v_{i})-\operatorname{div}%
_{\backslash e_{k}}g(v_{i})\right]}{2}\right)^{2}.
\end{align*}
Equality holds if and only if%
\begin{equation*}
u_{0}(v_{j})-\left(\operatorname{div}_{\backslash
e_{k}}g(v_{j})+g(e_{k})\right)=u_{0}(v_{i})-\left(\operatorname{div}_{\backslash
e_{k}}g(v_{i})-g(e_{k})\right),
\end{equation*}%
or equivalently
\begin{equation*}
g(e_{k})=\frac{\left[u_{0}(v_{j})-\operatorname{div}_{\backslash
e_{k}}g(v_{j})\right]-\left[u_{0}(v_{i})-\operatorname{div}_{\backslash e_{k}}g(v_{i})\right]}{2}%
=:Kg(e_{k}).
\end{equation*}%
Moreover, $\xi (x)$ is strictly convex and therefore strictly decreasing
for $x< Kg (e_{k})$ and strictly increasing for $x>Kg(e_{k})$%
. So the minimal value of $\xi (x)$ on the interval $\left[ -t,t\right] $
is only attained at
\begin{itemize}
	\item[(i)] the point $Kg(e_{k})$ if $Kg(e_{k})\in \left[ -t,t\right] $,
	\vskip1mm
	\item[(ii)] the point $-t$ if $Kg(e_{k})<-t$,
	\vskip1mm
	\item[(iii)] the point $t$ if $Kg(e_{k})>t$.
\end{itemize}
The assumption $\operatorname{div}g=\widetilde{\psi}$ then implies that 
$g(e_{k})$ must be the nearest point in the interval $\left[ -t,t\right] $
to $Kg(e_{k})$, implying that $T_{k}g(e_{k})=g(e_{k})$. Since $e_{k}\in E$ was
arbitrary, it follows that $T_{k}g\left( e_{k}\right) =g\left( e_{k}\right) $ for all
 $k=1,...,M$. Therefore $T_{k}g=g$ for all $k=1,...,M$ and we conclude that $Tg=g$.
\par Conversely, let us assume that $g\in t \mathcal{B}_{\ell^{\infty}\left(S_{E}\right)}$ and $Tg=g$. 
Then for any edge $e\in E$, $g(e)$ coincides with the point of the interval $\left[ -t,t\right] $ which is
nearest to $Kg(e)$. As $\left\Vert u_{0}-\operatorname{div}(\cdot)\right\Vert
_{\ell^{2}(S_{V})}$ is a convex function on $t\mathcal{B}_{\ell^{\infty}\left(S_{E}\right)}$, it is enough to show that $g$ minimizes $\left\Vert u_{0}-\operatorname{div}(\cdot)\right\Vert _{\ell^{2}(S_{V})}$ locally, i.e. it is enough to show that for some small $\varepsilon >0$ we have%
\begin{equation*}
\left\Vert u_{0}-\operatorname{div}g\right\Vert _{\ell^{2}(S_{V})}=\underset{\omega\in
D_{\varepsilon }}{\inf}\left\Vert u_{0}-\operatorname{div}\omega\right\Vert
_{\ell^{2}(S_{V})},
\end{equation*}%
where $D_{\varepsilon}$ is the tubular set given by $D_{\varepsilon }=\left\{ \omega\in t\mathcal{B}_{\ell^{\infty}\left(S_{E}\right)}:\ \left\Vert g-\omega\right\Vert _{\ell^{\infty }(S_{E})}\leq \varepsilon \right\}$.
Note that for any $\omega\in D_{\varepsilon }$ and $e\in E$ we have $\omega(e)\in \left[ -t,t\right] \cap \left[ g(e)-\varepsilon ,g(e)+\varepsilon \right]$.
The set $D_{\varepsilon }$ is a compact subset of $S_{E}$ and it therefore exists a function $\omega_{\varepsilon }\in D_{\varepsilon }$ such that%
\begin{equation} \label{eq:defomegaeps}
\left\Vert u_{0}-\operatorname{div}\omega_{\varepsilon }\right\Vert _{\ell^{2}(S_{V})}=%
\underset{\omega\in D_{\varepsilon }}{\inf }\left\Vert u_{0}-\operatorname{div}%
\omega\right\Vert _{\ell^{2}(S_{V})}.
\end{equation}%
So, we will need to prove that%
\begin{equation*}
\left\Vert u_{0}-\operatorname{div}g\right\Vert _{\ell^{2}(S_{V})}=\left\Vert u_{0}-%
\operatorname{div}\omega_{\varepsilon }\right\Vert _{\ell^{2}(S_{V})}.
\end{equation*}%
We first note that it follows from the necessity direction proved above that
for any edge $e\in E$, $\omega_{\varepsilon }(e)$ will coincide with the point of
the interval $\left[ -t,t\right] \cap \left[ g(e)-\varepsilon
,g(e)+\varepsilon \right] $, which is nearest to $K\omega_{\varepsilon }(e)$.

Let us now decompose the edge set $E$ into two parts. The first part denoted
 by $\Omega _{g}$ consists of the edges for which $Kg(e)$
does not belong to the interval $\left[-t,t\right] $, i.e. $\Omega_{g}=\left\{ e\in E:\ Kg(e)\notin \left[-t,t\right] \right\}$.
 As $g(e)$ is the nearest point in the interval $\left[ -t,t\right] $ to $%
Kg(e)$ we have
\begin{align*}
g(e)=\left\{
\begin{array}{rl}
-t & \text{if } Kg\left(e\right)< -t \\
+t & \text{if } Kg\left(e\right)> +t
\end{array} \right.\ \textrm{for edges}\ e\in \Omega _{u}.
\end{align*}
If the number $\varepsilon >0$ is small
enough, it follows from $\left\Vert g-\omega\right\Vert _{\ell^{\infty }(S_{E})}\leq
\varepsilon $ that on $e\in \Omega _{g}$ where we have $Kg(e)<-t$ we will
also have $K\omega_{\varepsilon }(e)<-t$ and therefore $\omega_{\varepsilon
}(e)=-t=g(e)$. Analogously, on $e\in
\Omega _{g}$ where $Kg(e)>t$ we will have $K\omega_{\varepsilon }(e)>t$ and
therefore $\omega_{\varepsilon }(e)=t=g(e)$. So we have
\begin{equation}
\label{12}
\omega_{\varepsilon }(e)=g(e)\ \textrm{for all}\ e\in \Omega _{g}.  
\end{equation}%
\par Next, we consider the remaining edges $E\backslash \Omega_{g}$. Let $G^{\prime}=(V,E\backslash \Omega
_{g})$, i.e. the graph $G$ with the edges in $\Omega _{g}$ removed. The graph $G^{\prime }$
is the union of several connected components $\left( V_{k},E_{k}\right),\ k=1,...,\ell$ 
 so that we have $V_{1}\cup ...\cup V_{\ell}=V$ and $E_{1}\cup ...\cup E_{\ell}=E\backslash \Omega _{g}$. 
Note that it is possible that some of the graphs $\left(V_{k},E_{k}\right)$ consist of just one single vertex. 
For these graphs there is nothing to prove because $E_{k}=\emptyset$. Let us now consider a subgraph $(V_{k},E_{k})$ where $E_{k}\neq \emptyset $.
On each $e\in E_{k}$ we have $Kg(e)\in \left[ -t,t\right] $ and therefore $%
g(e)=Kg(e)$, i.e. if $e=(v_{i},v_{j})$ then%
\begin{equation*}
g(e)=Kg(e)=\frac{\left[u_{0}(v_{j})-\operatorname{div}_{\backslash
e}u(v_{j})\right]-\left[u_{0}(v_{i})-\operatorname{div}_{\backslash e}u(v_{i})\right]}{2},
\end{equation*}%
or equivalently, in view of the definition of $\operatorname{div}_{\backslash e}g(\cdot)$, we get that
\begin{align*}
u_{0}(v_{j})-\operatorname{div}g(v_{j})=u_{0}(v_{i})-\operatorname{div}g(v_{i}).
\end{align*}
Note that operators $K$, $\operatorname{div}$ and $\operatorname{div%
}_{\backslash u(e)}$ are considered in the original setting of $G=(V,E)$. Therefore, for all $v\in V_{k}$ the values of $u_{0}(v)-\operatorname{div}g(v)$ are equal. It follows that
\begin{align}
\label{36}
\sum_{v\in V_{k}}\left[u_{0}(v)-\operatorname{div}g(v)\right]^{2}=
\left| V_{k}\right| \left( \frac{\sum_{v\in V_{k}}\left[ u_{0}(v)-%
\operatorname{div}g(v)\right]}{\left\vert V_{k}\right\vert }\right) ^{2}.  
\end{align}%
For $\omega_{\varepsilon }$ we can with Jensen's inequality derive the
corresponding inequality%
\begin{align}
\label{37}
\sum_{v\in V_{k}}\left[u_{0}(v)-\operatorname{div}\omega_{\varepsilon }(v)\right]^{2}\geq 
\left| V_{k}\right| \left( \frac{\sum_{v\in V_{k}}\left[u_{0}(v)-%
\operatorname{div}\omega_{\varepsilon }(v)\right]}{\left\vert V_{k}\right\vert }\right)^{2}. 
\end{align}%
Now, note that flows on edges in $E_{k}$ are canceled in the sums $\sum_{v\in V_{k}}\left[u_{0}(v)-\operatorname{div}g(v)\right]$ and $\sum_{v\in
V_{k}}\left[ u_{0}(v)-\operatorname{div}\omega_{\varepsilon }(v)\right]$. Therefore, only flows on edges in $\Omega _{g}$ remain in these sums. It then follows from \eqref{12} that 
\begin{equation*}
\sum_{v\in V_{k}}\left[u_{0}(v)-\operatorname{div}\omega_{\varepsilon }(v)\right]
=\sum_{v\in V_{k}}\left[u_{0}(v)-\operatorname{div}g(v)\right].
\end{equation*}
Therefore, taking into account \eqref{36} and \eqref{37}, we obtain
\begin{equation*}
\sum_{v\in V_{k}}\left[u_{0}(v)-\operatorname{div}\omega_{\varepsilon }(v)\right]^{2}\geq
\sum_{v\in V_{k}}\left[u_{0}(v)-\operatorname{div}g(v)\right]^{2}.
\end{equation*}
Summing over all $V_{k}$ gives
\begin{equation*}
\left\| u_{0}-\operatorname{div}\omega_{\varepsilon }\right\|_{\ell^{2}(S_{V})}^{2}\geq
\left\| u_{0}-\operatorname{div}g\right\|_{\ell^{2}(S_{V})}^{2},
\end{equation*}
and we conclude from the definition of $\omega_{\varepsilon }$, recall \eqref{eq:defomegaeps}, that
\begin{equation*}
\left\| u_{0}-\operatorname{div}\omega_{\varepsilon}\right\|_{\ell^{2}(S_{V})}^{2}=\left\| u_{0}-\operatorname{div}g\right\|_{\ell^{2}(S_{V})}^{2}.
\end{equation*}
So, $g$ minimizes $\left\Vert u_{0}-\operatorname{div}(\cdot)\right\Vert _{\ell^{2}(S_{V})}$ on $D_{\varepsilon }$ and therefore, by convexity, on $t\mathcal{B}_{\ell^{\infty }(S_{E})}$. Therefore, $Tg=g$ implies $\operatorname{div}g=\widetilde{\psi}$ and we have now established condition (1).
\par Finally, we prove condition (2). Note that by definition for $\forall g\in t\mathcal{B}_{\ell^{\infty }(S_{E})}$, the
operators $T_{k}$, $k=1,...,M$ satisfy
\begin{equation*}
\left\Vert u_{0}-\operatorname{div}\left(T_{k}u\right)\right\Vert _{\ell^{2}(S_{V})}\leq \left\Vert
u_{0}-\operatorname{div}u\right\Vert _{\ell^{2}(S_{V})},
\end{equation*}%
with equality if and only if $T_{k}g(e_{k})=g(e_{k})$. This implies that
\begin{equation*}
\left\Vert u_{0}-\operatorname{div}\left(Tg\right)\right\Vert _{\ell^{2}(S_{V})}\leq \left\Vert
u_{0}-\operatorname{div}g\right\Vert _{\ell^{2}(S_{V})},
\end{equation*}%
with equality if and only if $Tg=g$ which in turn by condition (1) is equivalent to $\operatorname{div}g=\widetilde{\psi}$. Hence for any $g\in t \mathcal{B}_{\ell^{\infty}\left(S_{E}\right)}$, $\textrm{if}\ \operatorname{div}g\neq\widetilde{\psi}$ then 
	\begin{align*}
		\left\| u_{0}-\operatorname{div}\left(Tg\right) \right\|_{\ell^{2}\left(S_{V}\right)} < \left\| u_{0}-\operatorname{div}g \right\|_{\ell^{2}\left(S_{V}\right)}.
	\end{align*}
\end{proof}
\vskip2mm
\par We are now ready to show the following theorem which establish that Algorithm \ref{ALGO} converges to the ROF-minimizer $u_{opt}$.
\vskip2mm
\begin{theorem}
\label{PropoTconv}
\par Let $g\in t \mathcal{B}_{\ell^{\infty}\left(S_{E}\right)}$ and $T:t \mathcal{B}_{\ell^{\infty}\left(S_{E}\right)} \rightarrow t \mathcal{B}_{\ell^{\infty}\left(S_{E}\right)}$ be the operator given by \eqref{U0}-\eqref{U1}. Then
\begin{align*}
	\operatorname{div}\left(T^{n}g\right) \rightarrow \widetilde{\psi}=u_{0}-u_{opt}\ \textrm{as}\ n\rightarrow +\infty.
\end{align*}
\end{theorem}
\vskip2mm
\begin{proof}
From Proposition \ref{TheOpTproperties}, it follows that $T$ is continuous and satisfies the conditions $(1)$ and $(2)$. 
These conditions in turn give that the sequence $\left(\left\|u_{0} - \operatorname{div}\left(T^{n}g\right) \right\|_{\ell^{2}(S_{V})}\right)_{n\in\mathbb{N}}$ is monotonically decreasing and bounded below by $\left\|u_{0}-\widetilde{\psi}\right\|_{\ell^{2}\left(S_{V}\right)}$. Therefore it converges. Let us now consider the sequence $\left(T^{n}g\right)_{n\in\mathbb{N}}\subset t \mathcal{B}_{\ell^{\infty}\left(S_{E}\right)}$. The ball $t \mathcal{B}_{\ell^{\infty}\left(S_{E}\right)}$ is a compact set and therefore has $\left(T^{n}g\right)_{n\in\mathbb{N}}$ a convergent subsequence in $t \mathcal{B}_{\ell^{\infty}\left(S_{E}\right)}$, say $\left(T^{n_{k}}g\right)_{k\in\mathbb{N}}$:
\begin{align*}
	\lim_{k\rightarrow \infty} T^{n_{k}}g = g_{\psi}\in t \mathcal{B}_{\ell^{\infty}\left(S_{E}\right)}.
\end{align*}
Since $T$, $\operatorname{div}$ and $\left\|\cdot\right\|_{\ell^{2}(S_V)}$ are continuous operators, we have
\begin{align*}	
  \left\| u_{0}-\operatorname{div}\left(Tg_{\psi}\right) \right\|_{\ell^{2}\left(S_{V}\right)}	&=\left\| u_{0}-\operatorname{div}\left(T \left(\lim_{k\rightarrow \infty} T^{n_{k}}g\right) \right) \right\|_{\ell^{2}\left(S_{V}\right)} \\
& =\lim_{k\rightarrow \infty}\left\| u_{0} - \operatorname{div}  \left(T \left(T^{n_{k}}g\right) \right) \right\|_{\ell^{2}\left(S_{V}\right)}\\&= \lim_{k\rightarrow \infty} \left\| u_{0} - \operatorname{div} \left(T^{n_{k}+1}g \right) \right\|_{\ell^{2}\left(S_{V}\right)}.
\end{align*}
As $T^{n_{k+1}}g=T^{m}T^{n_{k}+1}g$ for some $m\in\left\{0,1,2,...\right\}$, Proposition \ref{TheOpTproperties} implies
\begin{align*}
 \lim_{k\rightarrow \infty} \left\| u_{0} - \operatorname{div} \left(T^{n_{k}+1}g \right) \right\|_{\ell^{2}\left(S_{V}\right)}\geq \lim_{k\rightarrow \infty} \left\| u_{0} - \operatorname{div} \left(T^{n_{k+1}}g \right) \right\|_{\ell^{2}\left(S_{V}\right)}.
\end{align*}
The continuity of $\operatorname{div}$ and $\left\|\cdot\right\|_{\ell^{2}(S_V)}$ then gives
\begin{align*}	
  \left\| u_{0}-\operatorname{div}\left(Tg_{\psi}\right) \right\|_{\ell^{2}\left(S_{V}\right)}\geq\lim_{k\rightarrow \infty} \left\| u_{0} - \operatorname{div} \left(T^{n_{k+1}}g \right) \right\|_{\ell^{2}\left(S_{V}\right)}=\left\| u_{0} - \operatorname{div}g_{\psi} \right\|_{\ell^{2}\left(S_{V}\right)}.
\end{align*}
Applying Proposition \ref{TheOpTproperties} again, we conclude that
\begin{align*}
\operatorname{div}g_\psi = \widetilde{\psi}
\end{align*}
and therefore, by the continuity of $\operatorname{div}$,
\begin{align} \label{eq:limsubseq}
 \lim_{k\rightarrow\infty}\operatorname{div}\left(T^{n_{k}}g\right)=\widetilde{\psi}.
\end{align}
\par The final step is to show the convergence of the entire sequence $\left(\operatorname{div}\left(T^{n}g\right)\right)_{n\in\mathbb{N}}$. From \eqref{eq:limsubseq} follows that
\begin{align*}
	\lim_{k\rightarrow\infty}\left\| u_{0} - \operatorname{div}\left(T^{n_{k}}g\right) \right\|_{\ell^{2}(S_{V})}=\left\|u_{0}-\widetilde{\psi}\right\|_{\ell^{2}\left(S_{V}\right)}.
\end{align*}
Since the subsequence $\left(\left\| u_{0} - \operatorname{div}\left(T^{n_{k}}g\right)\right\|_{\ell^{2}(S_{V})}\right)_{k\in\mathbb{N}}$ must converge to the same limit as the convergent sequence $\left(\left\| u_{0} - \operatorname{div}\left(T^{n}g\right) \right\|_{\ell^{2}(S_{V})}\right)_{n\in\mathbb{N}}$, we conclude that
\begin{align*}
	\lim_{n\rightarrow\infty}\left\| u_{0} - \operatorname{div}\left(T^{n}g\right)\right\|_{\ell^{2}(S_{V})}=\left\|u_{0}-\widetilde{\psi}\right\|_{\ell^{2}\left(S_{V}\right)}.
\end{align*}
Therefore, as $\widetilde{\psi}$ is the unique nearest element to $u_{0}$ in $t\mathcal{B}_{BV^{\ast}(S_V)}=t \operatorname{div}\left(\mathcal{B}_{\ell^{\infty}\left(S_{E}\right)}\right)$, we have
\begin{align*}
	\lim_{n\rightarrow\infty}\operatorname{div}\left(T^{n}g\right)=\widetilde{\psi}.
\end{align*}
\end{proof}
\vskip2mm
\section{Numerical results}\label{SecNumExpts}
\par In order to evaluate the performance of the proposed method, we make a numerical comparison with other efficient iterative TV minimization methods and one of the best known state-of-the-art denoising methods, namely BM3D. We have tested different types of images and made comparisons in terms of convergence rate, peak signal to noise ratio (PSNR), running time and visual quality. All experiments were performed on a \textit{Windows 7 Professional 32-bit} computer with a \textit{Intel(R) Core(TM) i5-2400 CPU, 3.1 GHz} Processor and a RAM of \textit{4096 MB}.  
\subsection{Comparison with other iterative TV minimization methods}
In this subsection the proposed algorithm is compared numerically with two state-of-the-art iterative algorithms for TV denoising, the Split-Bregman algorithm \cite{TomGoldsteinStanleyOsher} and the Primal-Dual algorithm \cite{AChambollePock}. We include numerical results obtained by testing different types of images and various noise levels. More specifically we consider a denoising scenario of natural and cartoon images aimed to numerically evaluate and illustrate the proposed algorithm's convergence rate and PSNR. In our experiments we used Gaussian noise with standard deviation 10, 20 and 30.
\begin{itemize}
\item The implementation of the Split-Bregman algorithm was obtained from \cite{webSplitBregman}. To find the best performing regularization parameter $\lambda$ we performed a brute-force optimization in the interval $[1,30]$ uniformly quantized into 100 values. 
\vskip1mm
\item The implementation for the Primal-Dual algorithm was obtained from the publicly available repository GPU4Vision \url{https://github.com/VLOGroup/primal-dual-toolbox}. In this implementation $\tau=\sigma=1/\sqrt{8}, \gamma=0.7\lambda$, the value of $\lambda$ was optimized in the same range as the regularization parameter in the Split-Bregman algorithm. The parameter $\theta$ was dynamically updated at each iteration by the rule $\theta \leftarrow 1/\sqrt{1+2 \gamma \tau}$ as well as $\tau \leftarrow \tau \theta$ and $\sigma \leftarrow \sigma/\theta$.
\vskip1mm
\item We implemented the proposed ROF model on the graph Algorithm \ref{ALGO} and optimized the regularization parameter using the same parameter space as the Split-Bregman algorithm.
\end{itemize}
The stopping criteria for all approaches was set to $||u^k-u^{k-1}||/|| u^k || < 10^{-5}$, where $||\cdot||$ is the Frobenius norm and  $u^k$ is the current iterate of the numerical scheme.
\begin{figure}[t]
  \hfill
  \begin{tikzpicture}[zoomboxarray, zoomboxes below, zoomboxarray inner gap=0.1cm, zoomboxarray columns=1, zoomboxarray rows=1]
    \node [image node] { \includegraphics[width=0.30\textwidth]{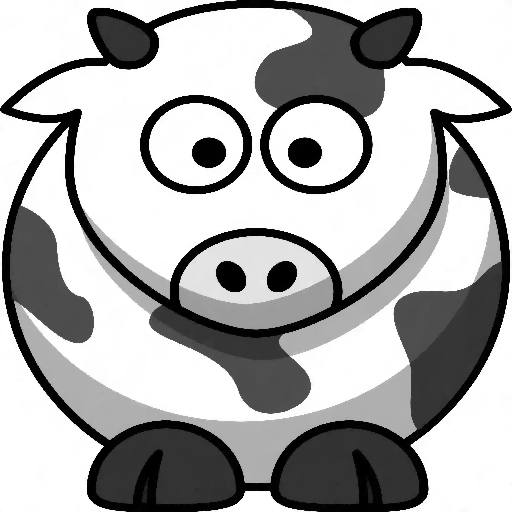} };
    \zoombox[color code=red,magnification=4]{0.43,0.45}
  \end{tikzpicture}
  \hfill \hfill
  \begin{tikzpicture}[zoomboxarray, zoomboxes below, zoomboxarray inner gap=0.1cm, zoomboxarray columns=1, zoomboxarray rows=1]
    \node [image node] { \includegraphics[width=0.30\textwidth]{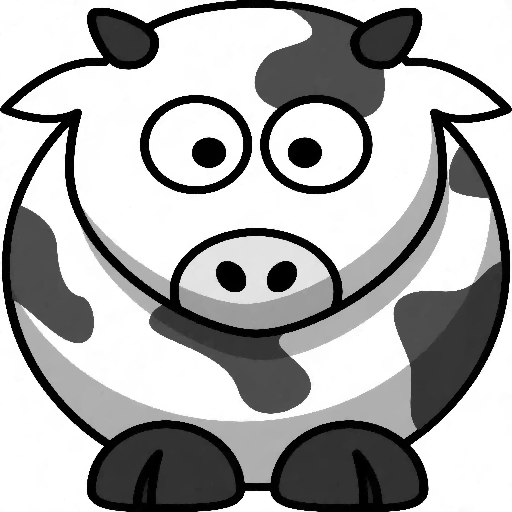} };
    \zoombox[color code=red,magnification=4]{0.43,0.45}
  \end{tikzpicture}
  \hfill \hfill
  \begin{tikzpicture}[zoomboxarray, zoomboxes below, zoomboxarray inner gap=0.1cm, zoomboxarray columns=1, zoomboxarray rows=1]
    \node [image node] { \includegraphics[width=0.30\textwidth]{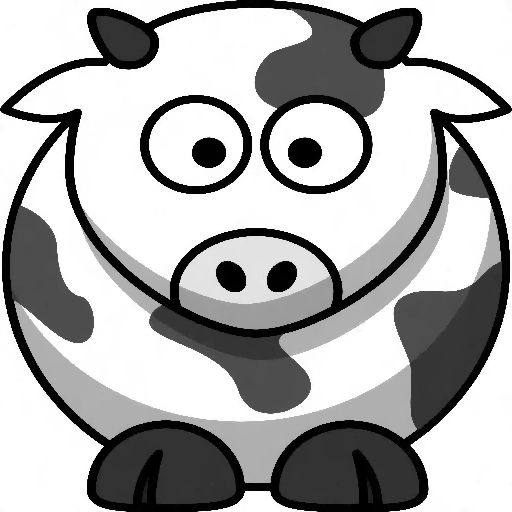} };
    \zoombox[color code=red,magnification=4]{0.43,0.45}
  \end{tikzpicture}
  \hfill \phantom{I} \\[3mm]
  \begin{tabularx}{1\linewidth}{XXX}
    (a) Split-Bregman (39.45)
    &
    (b) Primal-Dual (39.39)
    &
    (c) Proposed (40.80) \\
  \end{tabularx}
  \caption{Reconstructions when the original image was corrupted by Gaussian noise of standard deviation 10. Panels (a)-(c) show the images obtained at the best PSNR value for a brute-force parameter optimization strategy described in the main text. Obtained PSNR values are shown in parenthesis. Reaching the stopping criteria of $10^{-5}$, the proposed solution strategy shows improved PSNR compared to the other methods. We refer to Fig. \ref{fig:cartoon10} for empirical convergence results.}
  \label{fig:cartoon-images10}
\end{figure}

\begin{figure}[t]
  \centering
  \hfill \phantom{I}\\
  \includegraphics[width=0.49\linewidth]{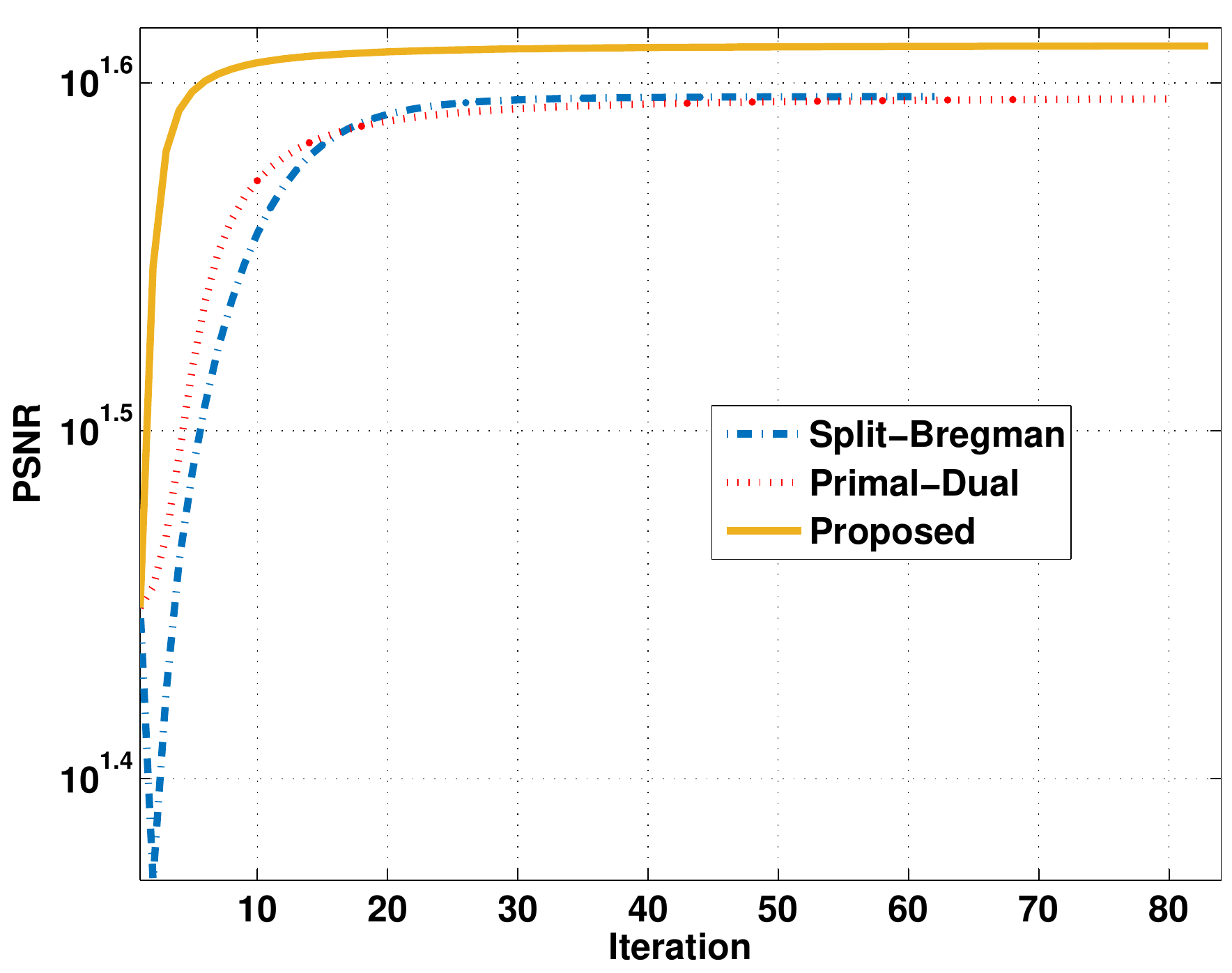} 
  \hfill\hfill
  \includegraphics[width=0.49\linewidth]{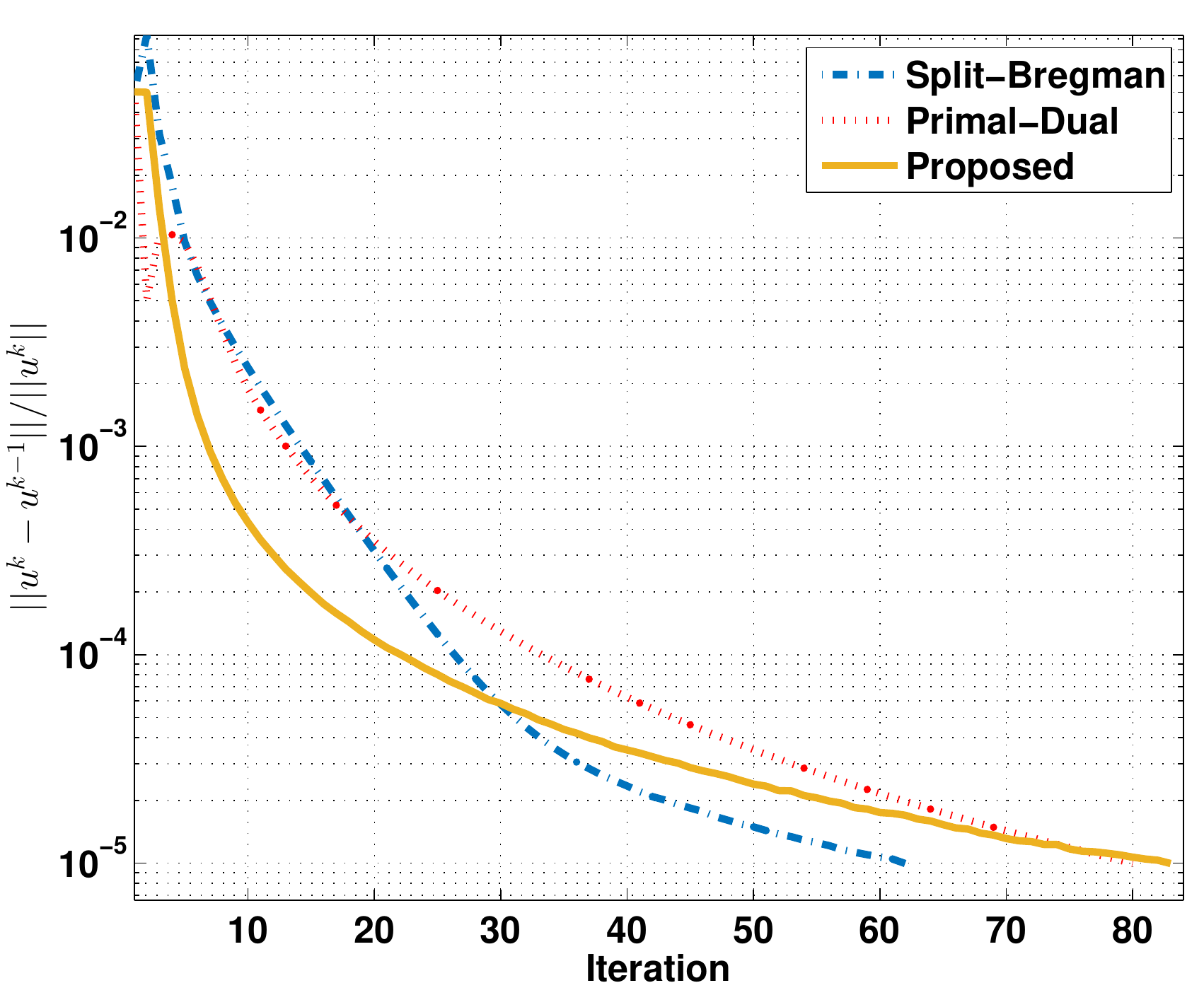} 
  \hfill\phantom{I}\\
  \begin{tabularx}{1\linewidth}{XX}
		(a) Empirical curves where PSNR is compared with iteration.
		&
    (b) Empirical convergence rate of the normalized error between two consecutive updates of the iterative schemes.
  \end{tabularx} 
\caption{Panels (a)-(b) show the PSNR and empirical convergence rates for the solution images in Fig. \ref{fig:cartoon-images10}. Panel (b) depicts that after 30 iterations the Split-Bregman scheme shows a smaller relative update between two consecutive update steps, this indicates earlier convergence for said method. However, after 30 iterations the improvement of the iterates is negligible compared to the ground truth data as seen in panel (a).}
\label{fig:cartoon10}
\end{figure}
\begin{figure}[t]
  \hfill
  \begin{tikzpicture}[zoomboxarray, zoomboxes below, zoomboxarray inner gap=0.1cm, zoomboxarray columns=1, zoomboxarray rows=1]
    \node [image node] { \includegraphics[width=0.30\textwidth]{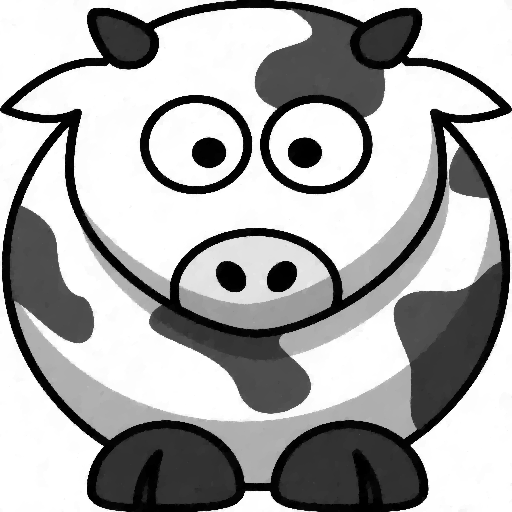} };
    \zoombox[color code=red,magnification=4]{0.43,0.45}
  \end{tikzpicture}
  \hfill \hfill
  \begin{tikzpicture}[zoomboxarray, zoomboxes below, zoomboxarray inner gap=0.1cm, zoomboxarray columns=1, zoomboxarray rows=1]
    \node [image node] { \includegraphics[width=0.30\textwidth]{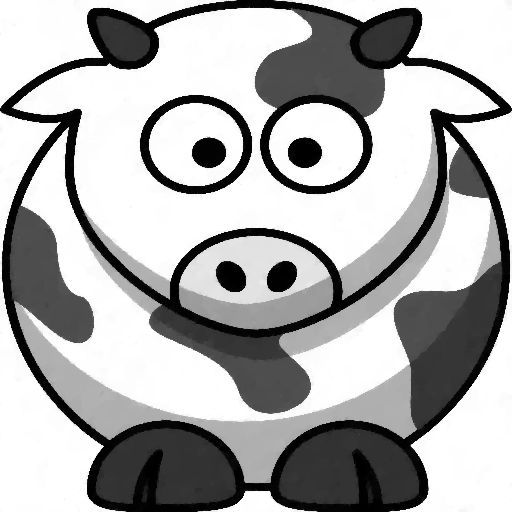} };
    \zoombox[color code=red,magnification=4]{0.43,0.45}
  \end{tikzpicture}
  \hfill \hfill
  \begin{tikzpicture}[zoomboxarray, zoomboxes below, zoomboxarray inner gap=0.1cm, zoomboxarray columns=1, zoomboxarray rows=1]
    \node [image node] { \includegraphics[width=0.30\textwidth]{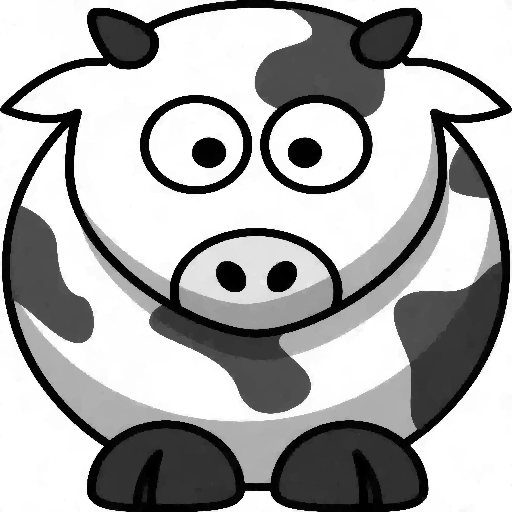} };
    \zoombox[color code=red,magnification=4]{0.43,0.45}
  \end{tikzpicture}
  \hfill \phantom{I} \\[3mm]
  \begin{tabularx}{1\linewidth}{XXX}
    (a) Split-Bregman (34.09)
    &
    (b) Primal-Dual (34.04)
    &
    (c) Proposed (35.17) \\
  \end{tabularx}
  \caption{Reconstructions when the original image was corrupted by Gaussian noise of standard deviation 20. Panels (a)-(c) show the images obtained at the best PSNR value for a brute-force parameter optimization strategy described in the main text. Obtained PSNR values are shown in parenthesis. Reaching the stopping criteria of $10^{-5}$, the proposed solution strategy shows improved PSNR compared to the other methods. We refer to Fig. \ref{fig:cartoon20} for empirical convergence results.}
  \label{fig:cartoon-images20}
\end{figure}
\begin{figure}[t]
  \centering
  \hfill \phantom{I}\\
  \includegraphics[width=0.49\linewidth]{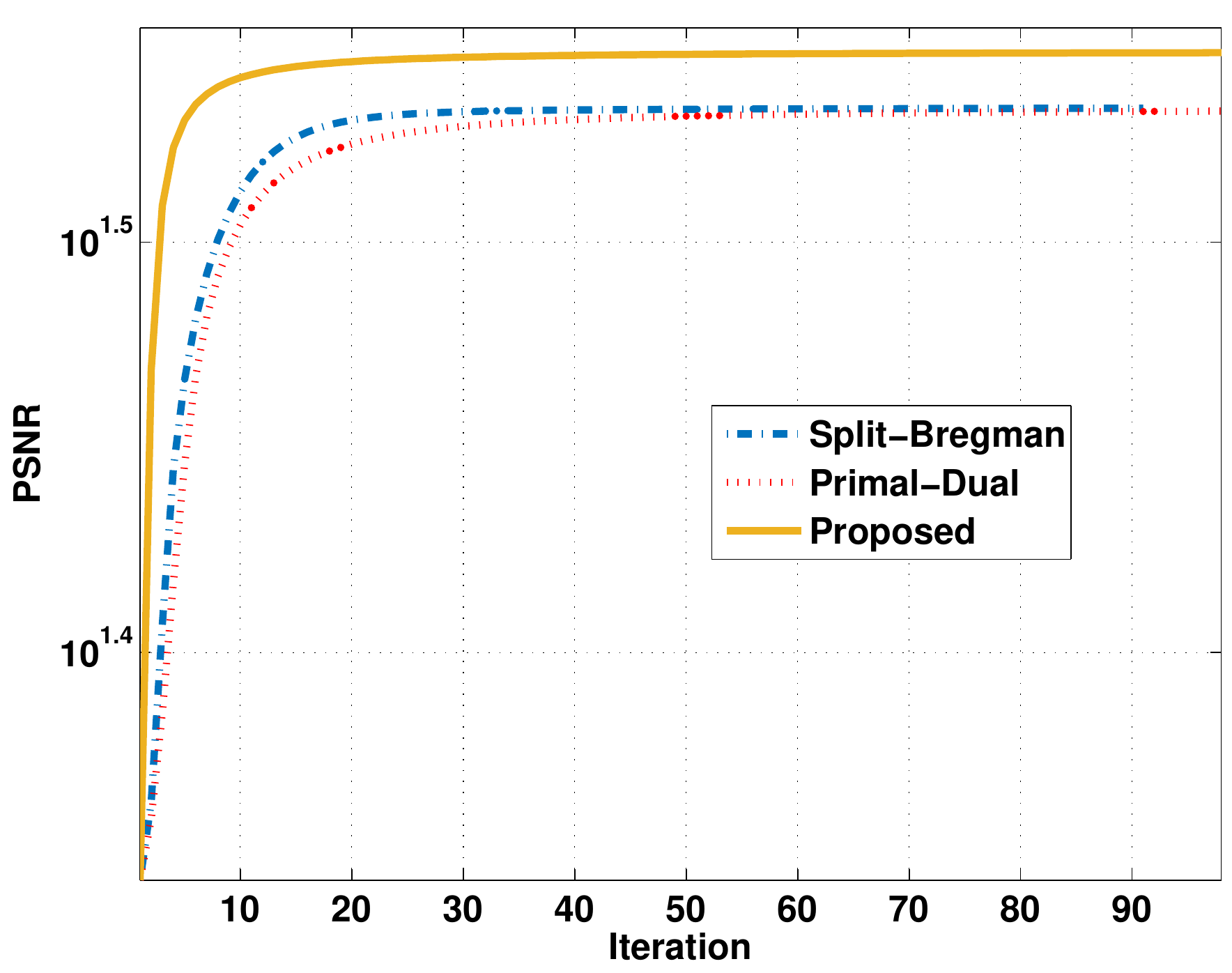} 
  \hfill\hfill
  \includegraphics[width=0.49\linewidth]{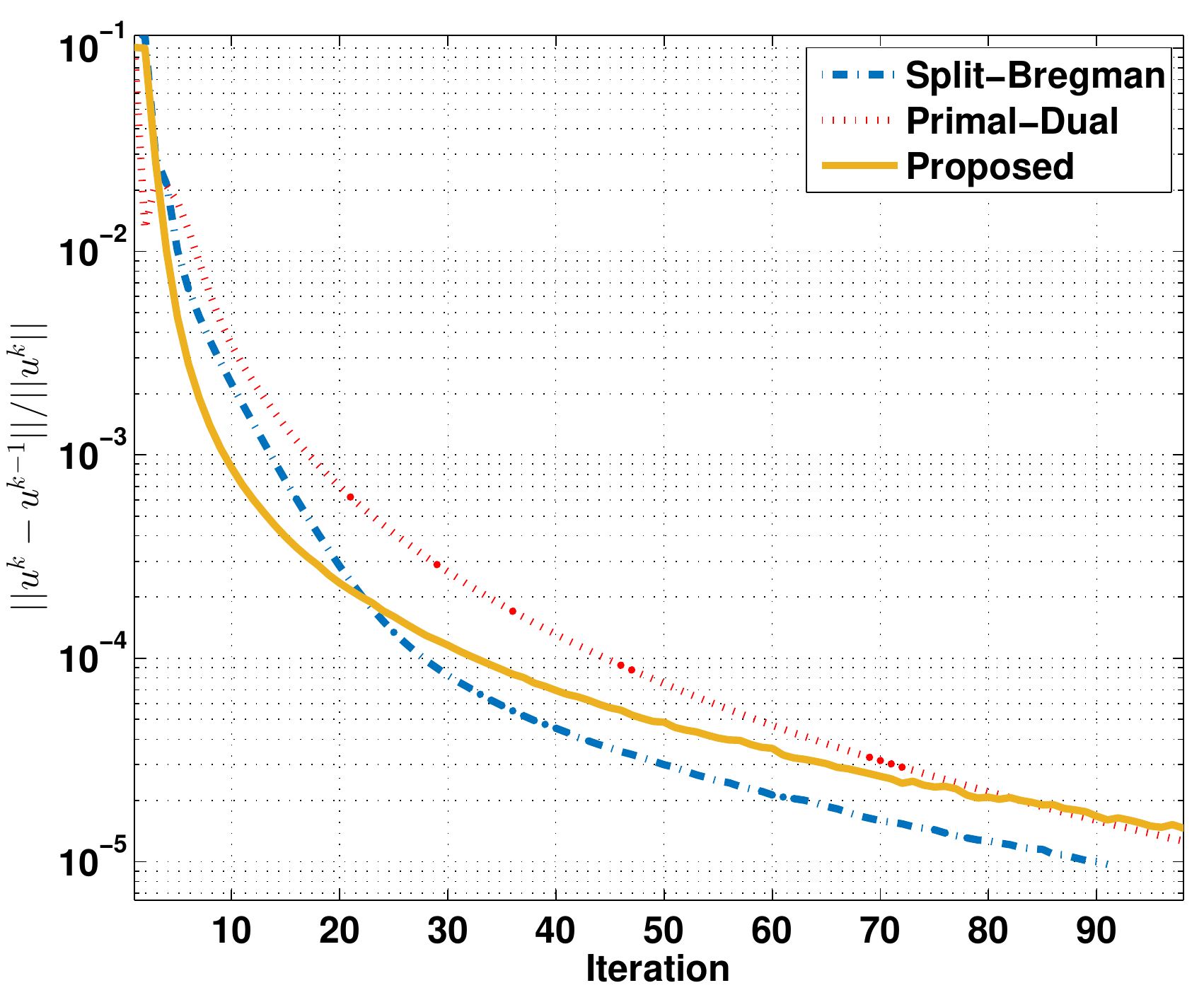} 
  \hfill\phantom{I}\\
  \begin{tabularx}{1\linewidth}{XX}
		(a) Empirical curves where PSNR is compared with iteration.
		&
    (b) Empirical convergence rate of the normalized error between two consecutive updates of the iterative schemes.  
  \end{tabularx} 
\caption{Panels (a)-(b) show the PSNR and empirical convergence rates for the solution images in Fig. \ref{fig:cartoon-images20}. Panel (b) depicts that after 20 iterations the Split-Bregman scheme shows a smaller relative update between two consecutive update steps, this indicates earlier convergence for said method. However, after 20 iterations the improvement of the iterates is negligible compared to the ground truth data as seen in panel (a).}
\label{fig:cartoon20}
\end{figure}
\begin{figure}[t]
  \hfill
  \begin{tikzpicture}[zoomboxarray, zoomboxes below, zoomboxarray inner gap=0.1cm, zoomboxarray columns=1, zoomboxarray rows=1]
    \node [image node] { \includegraphics[width=0.30\textwidth]{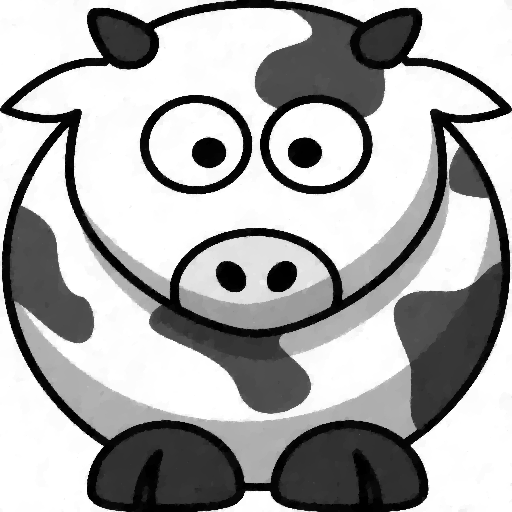} };
    \zoombox[color code=red,magnification=4]{0.43,0.45}
  \end{tikzpicture}
  \hfill \hfill
  \begin{tikzpicture}[zoomboxarray, zoomboxes below, zoomboxarray inner gap=0.1cm, zoomboxarray columns=1, zoomboxarray rows=1]
    \node [image node] { \includegraphics[width=0.30\textwidth]{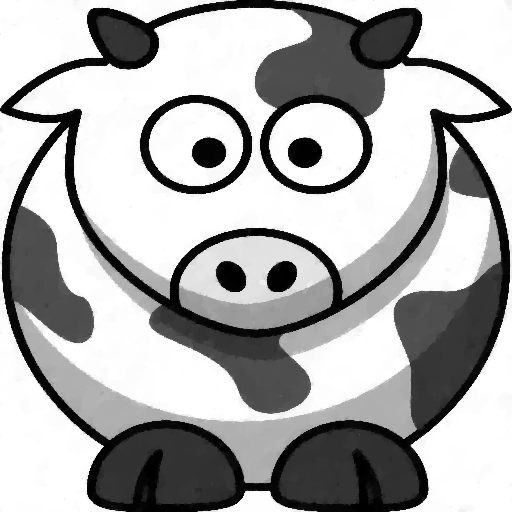} };
    \zoombox[color code=red,magnification=4]{0.43,0.45}
  \end{tikzpicture}
  \hfill \hfill
  \begin{tikzpicture}[zoomboxarray, zoomboxes below, zoomboxarray inner gap=0.1cm, zoomboxarray columns=1, zoomboxarray rows=1]
    \node [image node] { \includegraphics[width=0.30\textwidth]{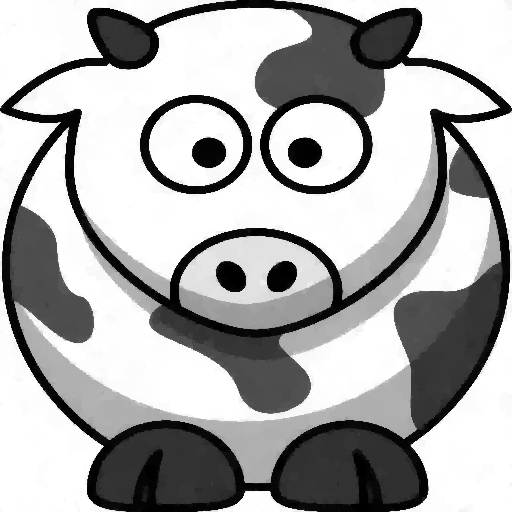} };
    \zoombox[color code=red,magnification=4]{0.43,0.45}
  \end{tikzpicture}
  \hfill \phantom{I} \\[3mm]
  \begin{tabularx}{1\linewidth}{XXX}
    (a) Split-Bregman (31.09)
    &
    (b) Primal-Dual (31.05)
    &
    (c) Proposed (31.94) \\
  \end{tabularx}
  \caption{Reconstructions when the original image was corrupted by Gaussian noise of standard deviation 30. Panels (a)-(c) show the images obtained at the best PSNR value for a brute-force parameter optimization strategy described in the main text. Obtained PSNR values are shown in parenthesis. Reaching the stopping criteria of $10^{-5}$, the proposed solution strategy shows improved PSNR compared to the other methods. We refer to Fig. \ref{fig:cartoon30} for empirical convergence results.}
  \label{fig:cartoon-images30}
\end{figure}

\begin{figure}[t]
  \centering
  \hfill \phantom{I}\\
  \includegraphics[width=0.49\linewidth]{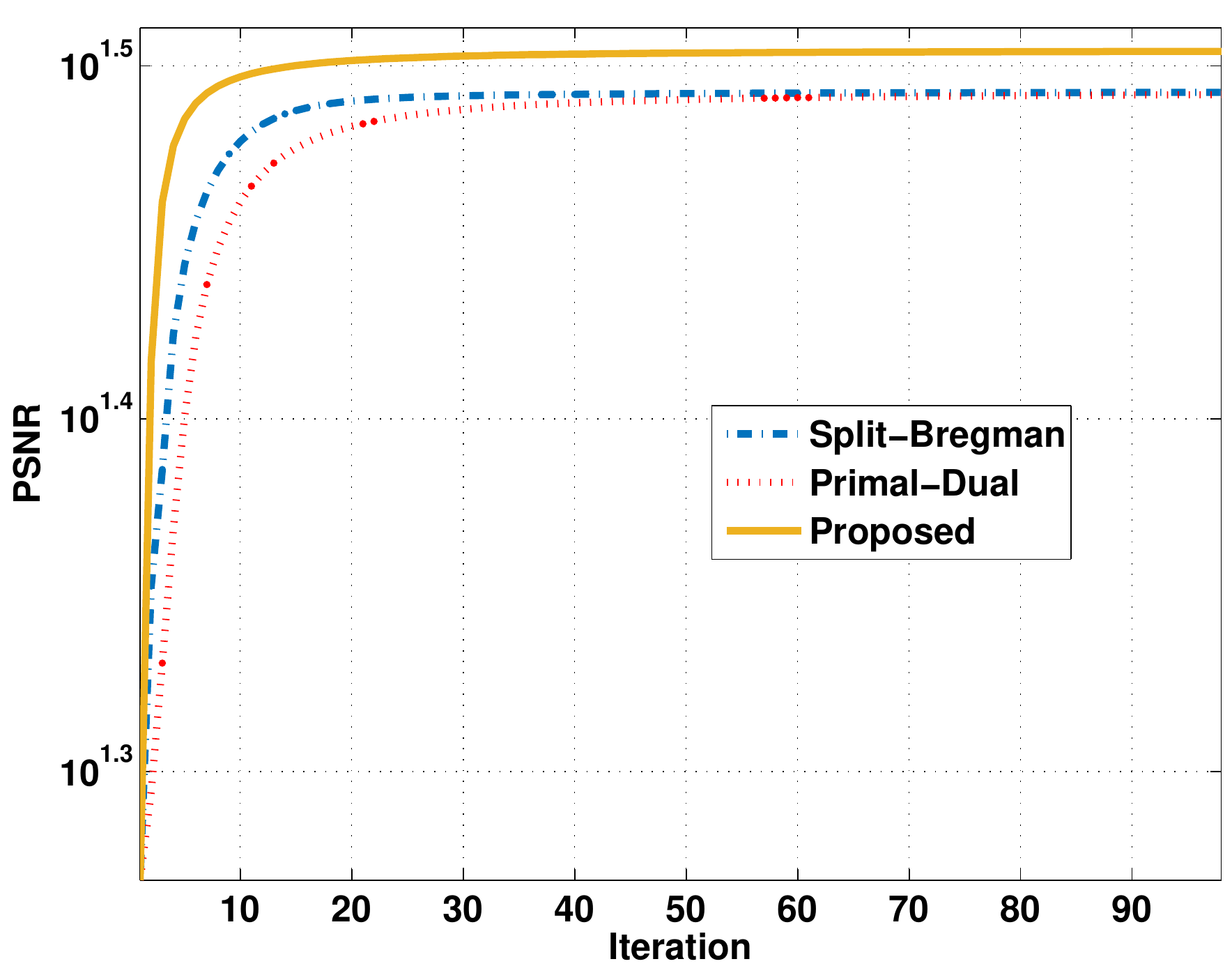} 
  \hfill\hfill
  \includegraphics[width=0.49\linewidth]{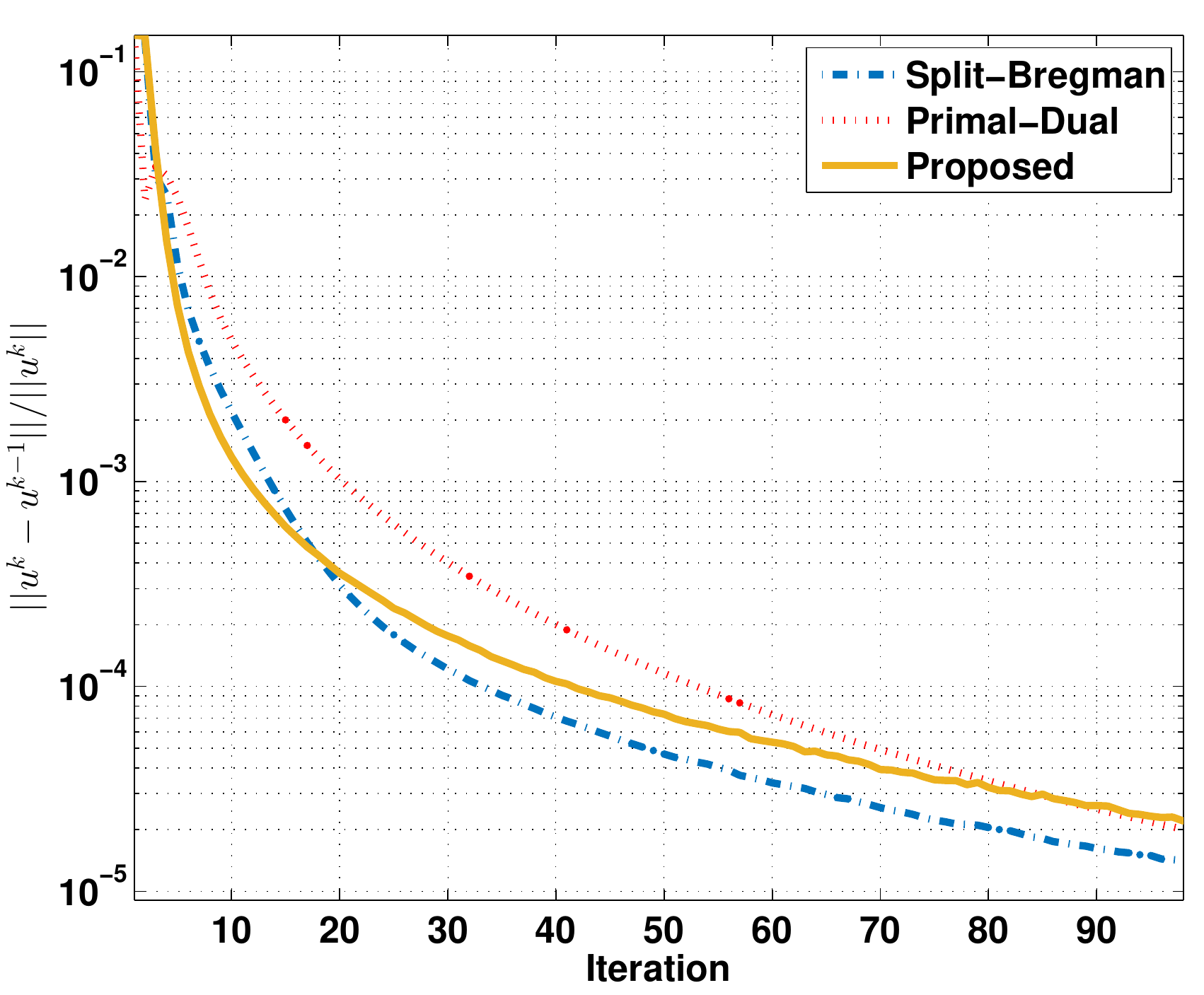} 
  \hfill\phantom{I}\\
  \begin{tabularx}{1\linewidth}{XX}
		(a) Empirical curves where PSNR is compared with iteration.
		&
    (b) Empirical convergence rate of the normalized error between two consecutive updates of the iterative schemes. 
  \end{tabularx} 
\caption{Panels (a)-(b) show the PSNR and empirical convergence rates for the solution images in Fig. \ref{fig:cartoon-images30}. Panel (b) depicts that after 20 iterations the Split-Bregman scheme shows a smaller relative update between two consecutive update steps, this indicates earlier convergence for said method. However, after 20 iterations the improvement of the iterates is negligible compared to the ground truth data as seen in panel (a). }
\label{fig:cartoon30}
\end{figure}
\subsection*{Cartoon image denoising}
In this example we have the \emph{exact} ground truth image data $u$ available which makes an objective evaluation of the methods possible. Figures \ref{fig:cartoon-images10}, \ref{fig:cartoon-images20} and \ref{fig:cartoon-images30} show the qualitative results for each evaluated method and noise level with the corresponding best PSNR values. The proposed method produces results with the best PSNR value in all cases. The visual quality of the results produced by all methods is comparable. Panels (a) of Figures \ref{fig:cartoon10}, \ref{fig:cartoon20} and \ref{fig:cartoon30} show the respective peak signal to noise ratio (PSNR) curves for the highest obtained PSNR values, obtained after a dense parameter grid search as previously described, for the respective methods and noise levels. Panels (b) of the same figures show the descent towards the stopping criteria. Each algorithm was terminated when the normalized difference between the current iterate and the previous iterate became smaller than $10^{-5}$, this is illustrated in panels (b). From the same figures it is clear that up till 30 (resp. 20) iterations for noise level 10 (resp. for noise levels 20 and 30) the proposed method shows a faster convergence rate. Note, however, that after these number of iterations any further updates of the iterative schemes have an negligible effect to the end result.
\begin{figure}[t]
  \hfill
  \begin{tikzpicture}[zoomboxarray, zoomboxes below, zoomboxarray inner gap=0.1cm, zoomboxarray columns=1, zoomboxarray rows=1]
    \node [image node] { \includegraphics[width=0.30\textwidth]{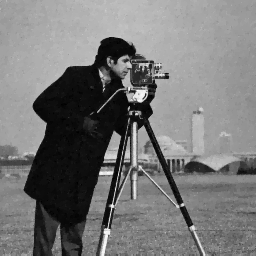} };
    \zoombox[color code=red,magnification=4]{0.48,0.75}
  \end{tikzpicture}
  \hfill \hfill
  \begin{tikzpicture}[zoomboxarray, zoomboxes below, zoomboxarray inner gap=0.1cm, zoomboxarray columns=1, zoomboxarray rows=1]
    \node [image node] { \includegraphics[width=0.30\textwidth]{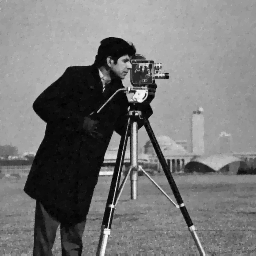} };
    \zoombox[color code=red,magnification=4]{0.48,0.75}
  \end{tikzpicture}
  \hfill \hfill
  \begin{tikzpicture}[zoomboxarray, zoomboxes below, zoomboxarray inner gap=0.1cm, zoomboxarray columns=1, zoomboxarray rows=1]
    \node [image node] { \includegraphics[width=0.30\textwidth]{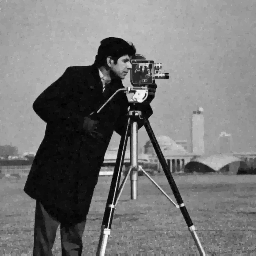} };
    \zoombox[color code=red,magnification=4]{0.48,0.75}
  \end{tikzpicture}
  \hfill \phantom{I} \\[3mm]
  \begin{tabularx}{1\linewidth}{XXX}
    (a) Split-Bregman (32.73)
    &
    (b) Primal-Dual (32.72)
    &
    (c) Proposed (32.73) \\
  \end{tabularx}
  \caption{Reconstructions when the original image was corrupted by Gaussian noise of standard deviation 10. Panels (a)-(c) show the images obtained at the best PSNR value for a brute-force parameter optimization strategy described in the main text. Obtained PSNR values are shown in parenthesis. The methods produces near identical PSNR values. We refer to Fig. \ref{fig:natural10} for empirical convergence results.}
  \label{fig:natural-images10}
\end{figure}
\begin{figure}[t]
  \centering
  \hfill \phantom{I}\\
  \hfill
  \includegraphics[width=0.49\linewidth]{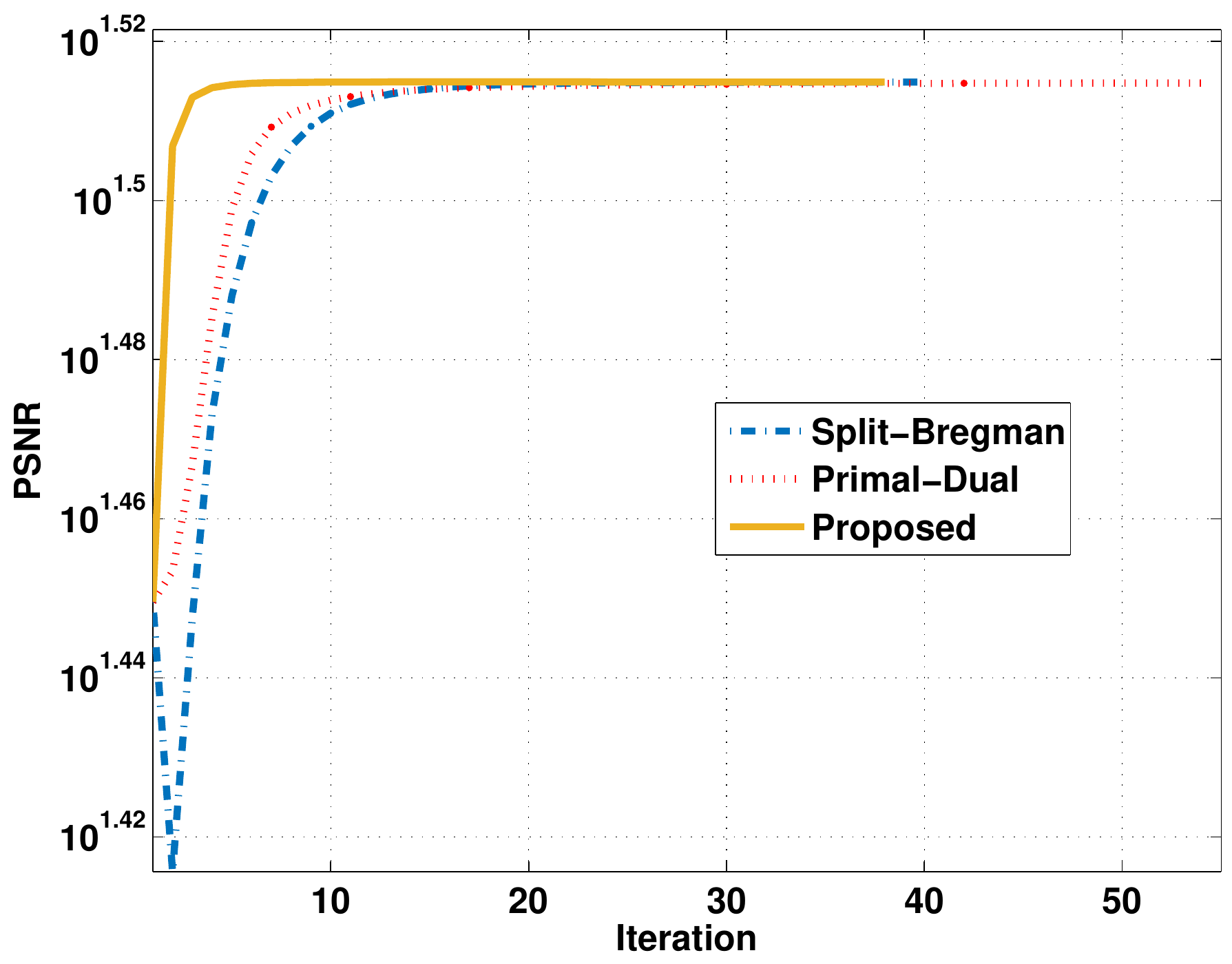} 
  \hfill\hfill
  \includegraphics[width=0.49\linewidth]{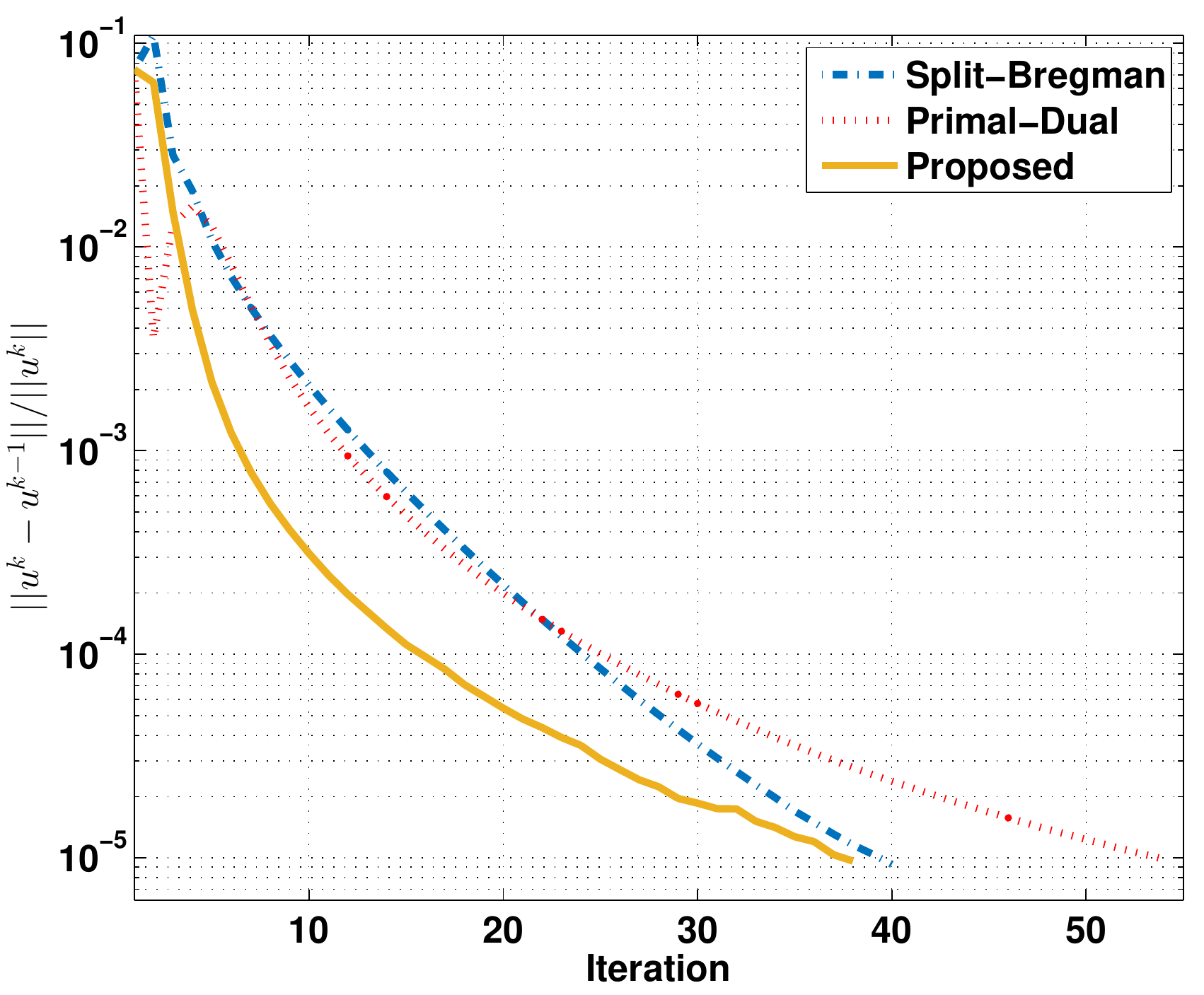} 
  \hfill \phantom{I}\\
  \begin{tabularx}{1\linewidth}{XXX}
			(a) Empirical curves where PSNR is compared with iteration.
			&
      (b) Empirical convergence rate of the normalized error between two consecutive updates of the iterative schemes.
  \end{tabularx}
    \caption{Panels (a)-(b) show the PSNR and empirical convergence rates to obtain the solution images in Fig. \ref{fig:natural-images10}. In this example, the proposed solution scheme shows an improved convergence rate compared to the Split-Bregman and the Primal-Dual approaches, yet resulting in near identical error values as seen in panel (a). }
  \label{fig:natural10}
\end{figure}
\begin{figure}[t]
  \hfill
  \begin{tikzpicture}[zoomboxarray, zoomboxes below, zoomboxarray inner gap=0.1cm, zoomboxarray columns=1, zoomboxarray rows=1]
    \node [image node] { \includegraphics[width=0.30\textwidth]{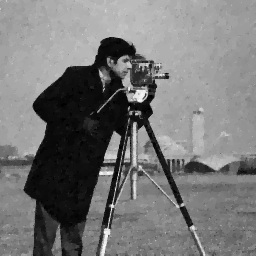} };
    \zoombox[color code=red,magnification=4]{0.48,0.75}
  \end{tikzpicture}
  \hfill \hfill
  \begin{tikzpicture}[zoomboxarray, zoomboxes below, zoomboxarray inner gap=0.1cm, zoomboxarray columns=1, zoomboxarray rows=1]
    \node [image node] { \includegraphics[width=0.30\textwidth]{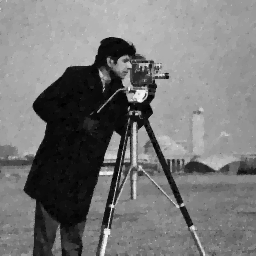} };
    \zoombox[color code=red,magnification=4]{0.48,0.75}
  \end{tikzpicture}
  \hfill \hfill
  \begin{tikzpicture}[zoomboxarray, zoomboxes below, zoomboxarray inner gap=0.1cm, zoomboxarray columns=1, zoomboxarray rows=1]
    \node [image node] { \includegraphics[width=0.30\textwidth]{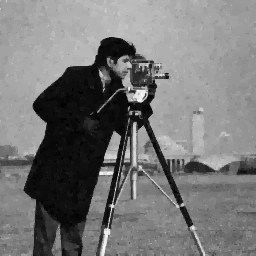} };
    \zoombox[color code=red,magnification=4]{0.48,0.75}
  \end{tikzpicture}
  \hfill \phantom{I} \\[3mm]
  \begin{tabularx}{1\linewidth}{XXX}
    (a) Split-Bregman (28.93)
    &
    (b) Primal-Dual (28.92)
    &
    (c) Proposed (28.93) \\
  \end{tabularx}
  \caption{Reconstructions when the original image was corrupted by Gaussian noise of standard deviation 20. Panels (a)-(c) show the images obtained at the best PSNR value for a brute-force parameter optimization strategy described in the main text. Obtained PSNR values are shown in parenthesis. The methods produces near identical PSNR values. We refer to Fig. \ref{fig:natural20} for empirical convergence results.}
  \label{fig:natural-images20}
\end{figure} 
\begin{figure}[t]
  \centering
  \hfill \phantom{I}\\
  \hfill
  \includegraphics[width=0.49\linewidth]{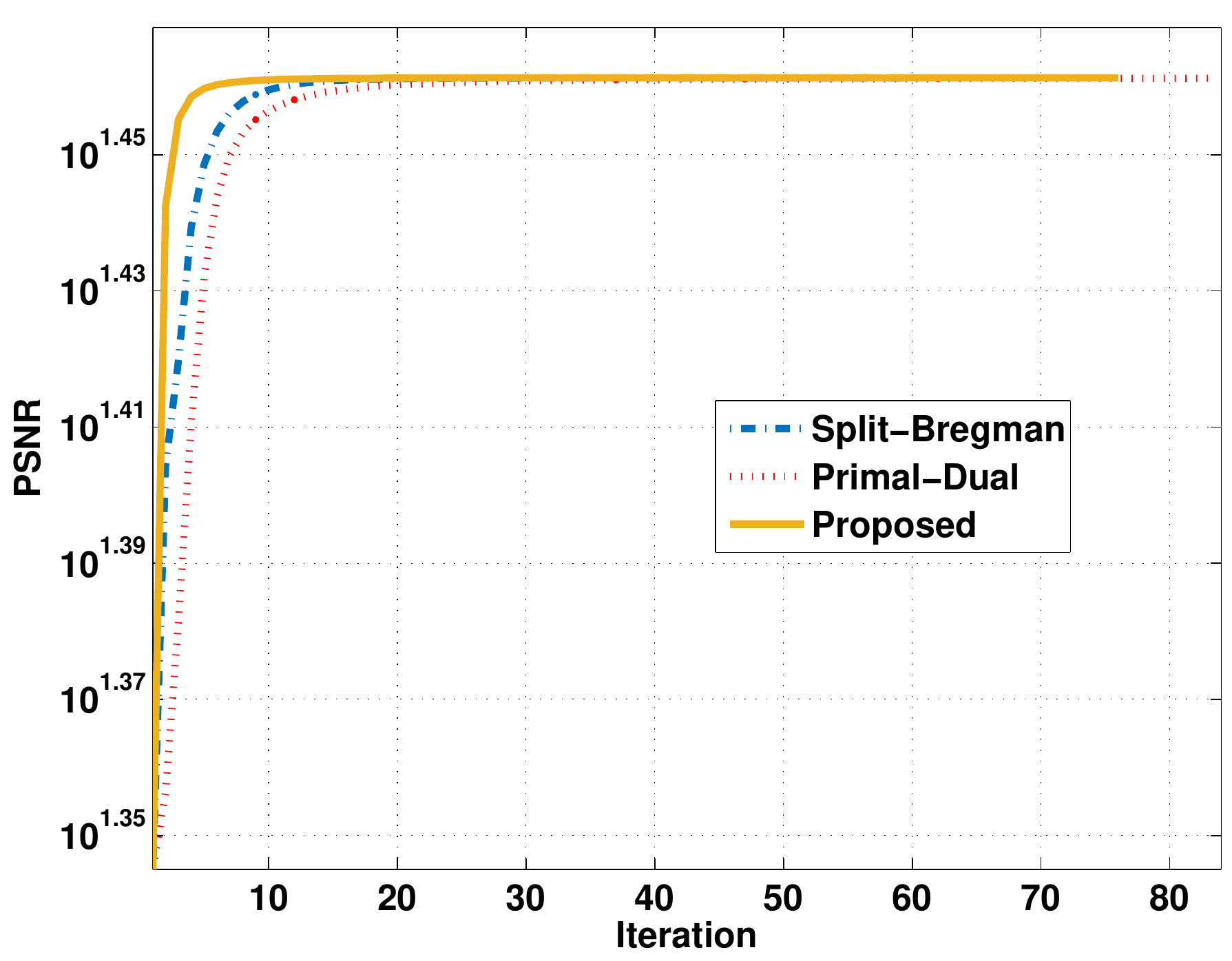} 
  \hfill\hfill
  \includegraphics[width=0.49\linewidth]{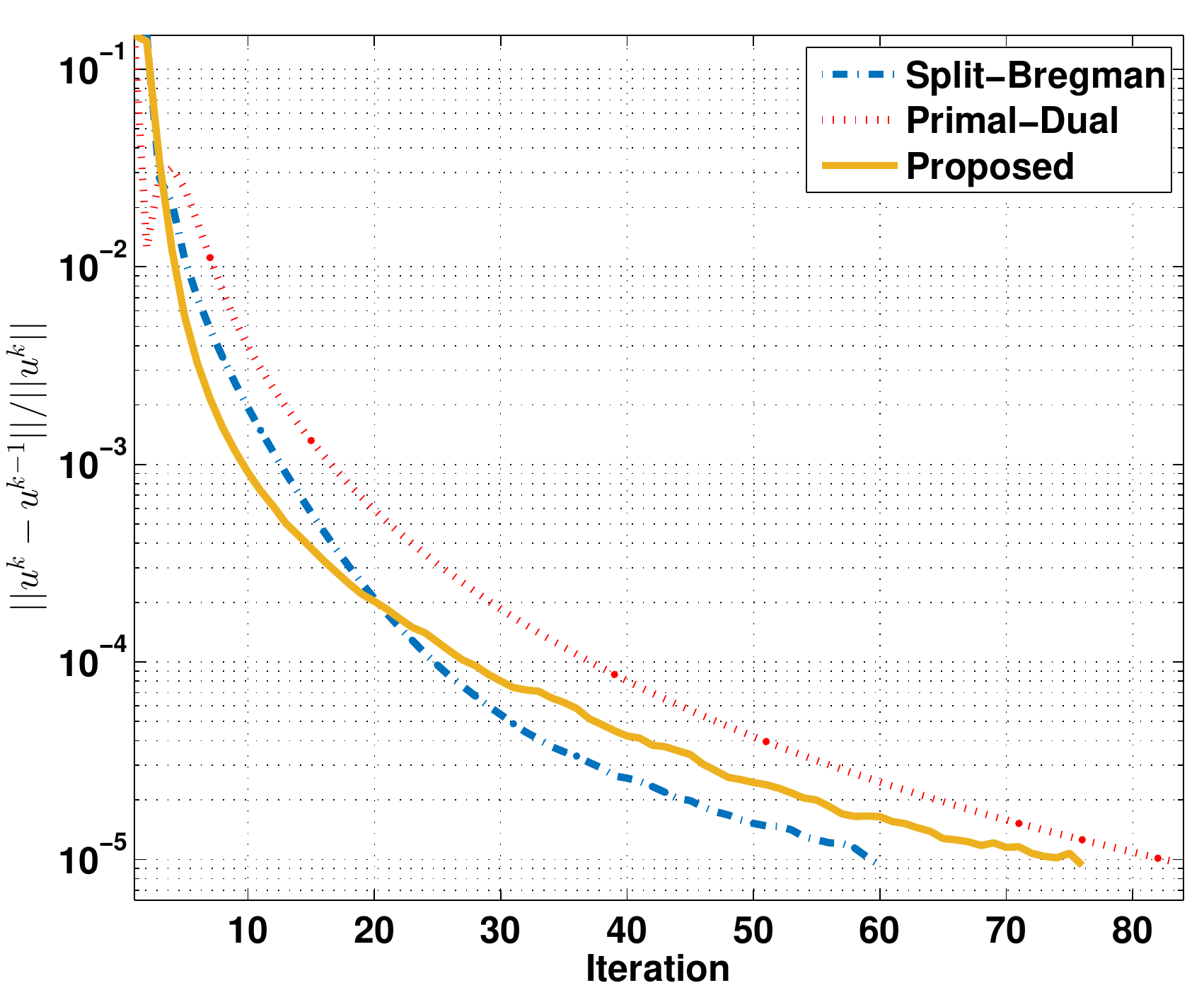} 
  \hfill \phantom{I}\\
  \begin{tabularx}{1\linewidth}{XXX}
			(a) Empirical curves where PSNR is compared with iteration.
			&
      (b) Empirical convergence rate of the normalized error between two consecutive updates of the iterative schemes.
  \end{tabularx}
    \caption{Panels (a)-(b) show the PSNR and empirical convergence rates to obtain the solution images in Fig. \ref{fig:natural-images20}. In this example, the proposed solution scheme shows an improved convergence rate compared to the Split-Bregman and the Primal-Dual approaches, yet resulting in near identical error values as seen in panel(a). }
  \label{fig:natural20}
\end{figure}
\begin{figure}[t]
  \hfill
  \begin{tikzpicture}[zoomboxarray, zoomboxes below, zoomboxarray inner gap=0.1cm, zoomboxarray columns=1, zoomboxarray rows=1]
    \node [image node] { \includegraphics[width=0.30\textwidth]{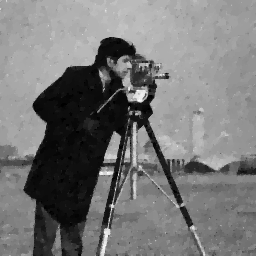} };
    \zoombox[color code=red,magnification=4]{0.48,0.75}
  \end{tikzpicture}
  \hfill \hfill
  \begin{tikzpicture}[zoomboxarray, zoomboxes below, zoomboxarray inner gap=0.1cm, zoomboxarray columns=1, zoomboxarray rows=1]
    \node [image node] { \includegraphics[width=0.30\textwidth]{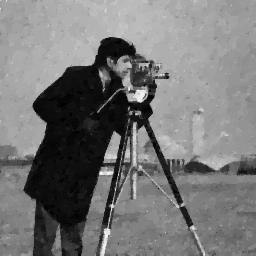} };
    \zoombox[color code=red,magnification=4]{0.48,0.75}
  \end{tikzpicture}
  \hfill \hfill
  \begin{tikzpicture}[zoomboxarray, zoomboxes below, zoomboxarray inner gap=0.1cm, zoomboxarray columns=1, zoomboxarray rows=1]
    \node [image node] { \includegraphics[width=0.30\textwidth]{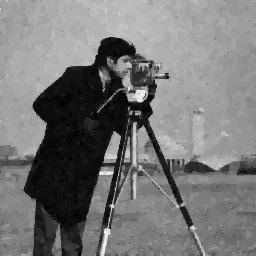} };
    \zoombox[color code=red,magnification=4]{0.48,0.75}
  \end{tikzpicture}
  \hfill \phantom{I} \\[3mm]
  \begin{tabularx}{1\linewidth}{XXX}
    (a) Split-Bregman (26.95)
    &
    (b) Primal-Dual (26.95)
    &
    (c) Proposed (26.94) \\
  \end{tabularx}
  \caption{Reconstructions when the original image was corrupted by Gaussian noise of standard deviation 30. Panels (a)-(c) show the images obtained at the best PSNR value for a brute-force parameter optimization strategy described in the main text. Obtained PSNR values are shown in parenthesis. The methods produces near identical PSNR values. We refer to Fig. \ref{fig:natural30} for empirical convergence results.}
  \label{fig:natural-images30}
\end{figure}
\begin{figure}[t]
  \centering
  \hfill \phantom{I}\\
  \hfill
  \includegraphics[width=0.49\linewidth]{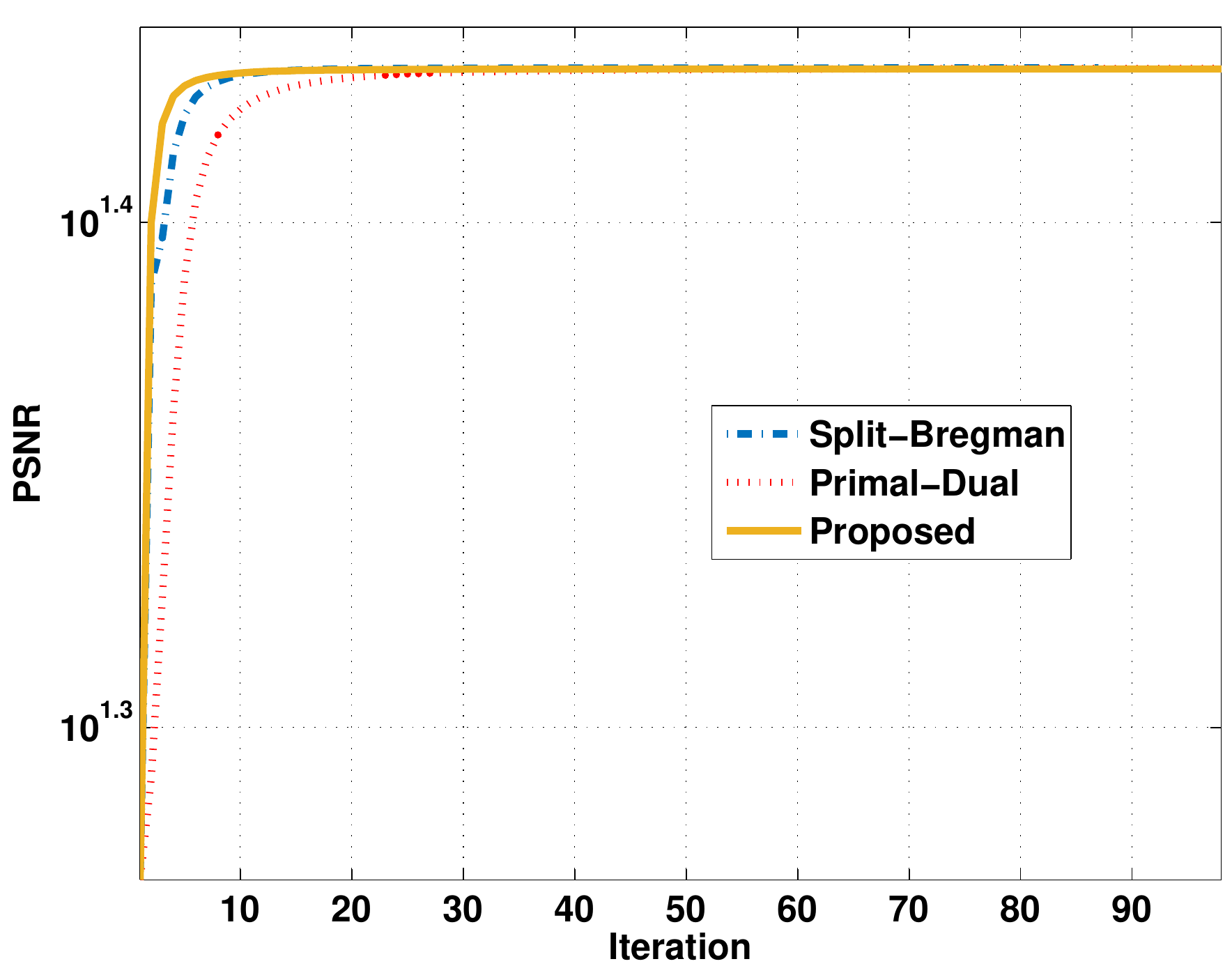} 
  \hfill\hfill
  \includegraphics[width=0.49\linewidth]{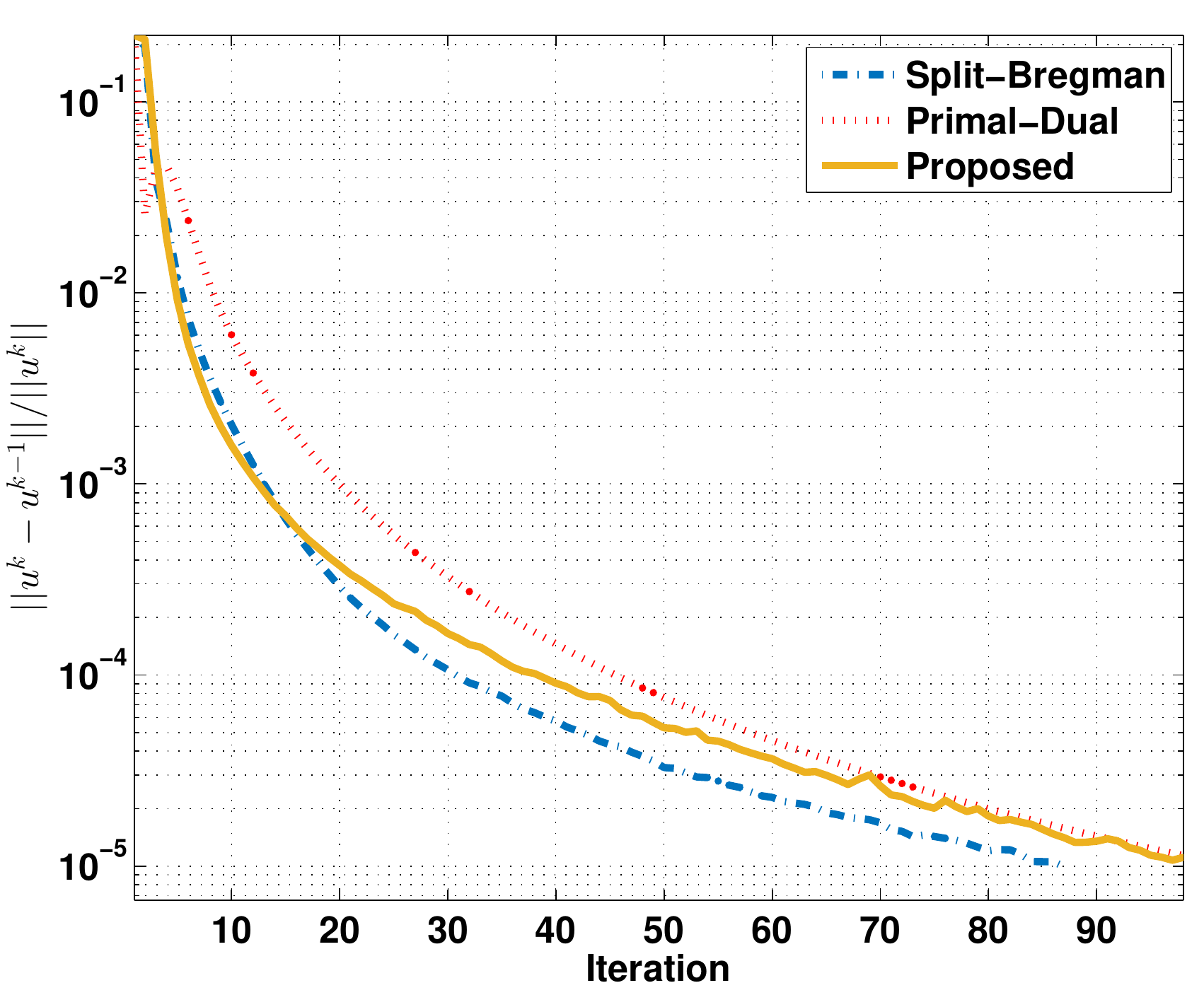} 
  \hfill \phantom{I}\\
  \begin{tabularx}{1\linewidth}{XXX}
			(a) Empirical curves where PSNR is compared with iteration.
			&
      (b) Empirical convergence rate of the normalized error between two consecutive updates of the iterative schemes. 
  \end{tabularx}
  \caption{Panels (a)-(b) show the PSNR and empirical convergence rates to obtain the solution images in Fig. \ref{fig:natural-images30}. In this example, the proposed solution scheme shows an improved convergence rate compared to the Split-Bregman and the Primal-Dual approaches, yet resulting in near identical error values as seen in panel (a).}
  \label{fig:natural30}
\end{figure}
\subsection*{Natural image denoising}
In this imaging scenario, we denoise the ``cameraman'' image. We have the \emph{exact} ground truth image data $u$ available so that it is possible to evaluate the methods objectively. Figures \ref{fig:natural-images10}, \ref{fig:natural-images20} and \ref{fig:natural-images30} show the qualitative results for each evaluated method and noise level with the corresponding best PSNR values. The proposed method produces results with the best or comparable PSNR value in all cases. The visual quality of the results produced by all methods is also comparable. In accordance with the previous example, panels (a) of Figures \ref{fig:natural10}, \ref{fig:natural20} and \ref{fig:natural30} show the respective peak signal to noise ratio (PSNR) curves for the highest obtained PSNR values, obtained after a dense parameter grid search as previously described, for the respective methods and noise levels. Further, panels (b) of the same figures show the descent towards the stopping criteria. Each algorithm was terminated when the normalized difference between the current iterate and the previous iterate became smaller than $10^{-5}$, this is illustrated in panels (b). From the same figures it is clear that similarly to the cartoon image, the relative convergence rate seen in (b) is initially faster for the proposed method than for the compared methods, and then it is overcome by the Split-Bregman algorithm at a certain number of iterations after which the improvement of the image quality is insignificant for all methods. Panels (b) shows that the best error rates for the proposed method are obtained at comparable fewer iterations than the Split-Bregman and the Primal-Dual approaches.  
\subsection{Comparison with BM3D}
\begin{figure}[t]
\centering
\begin{subfigure}{.475\linewidth}
  \centering
  \includegraphics[width=\linewidth]{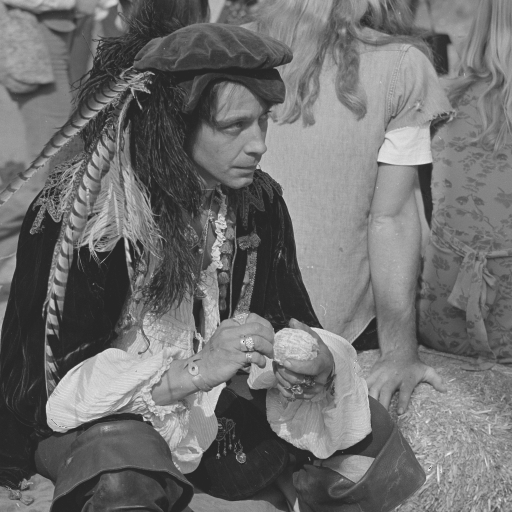}
  \caption{Man}
  \label{sfig:Man}
\end{subfigure}%
\hspace{1em}%
\begin{subfigure}{.475\linewidth}
  \centering
  \includegraphics[width=\linewidth]{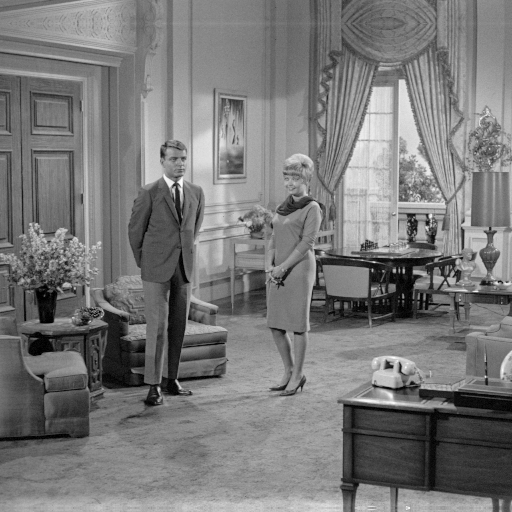}
  \caption{Couple}
  \label{sfig:Couple}
\end{subfigure}\\
\vspace{1em}
\begin{subfigure}{.475\linewidth}
  \centering
  \includegraphics[width=\linewidth]{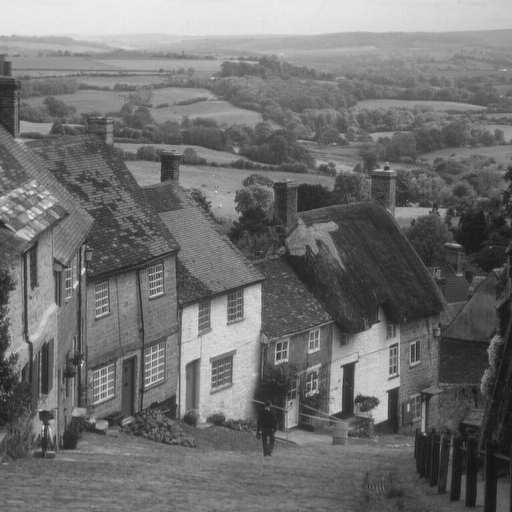}
  \caption{Hill}
  \label{sfig:Hill}
\end{subfigure}%
\hspace{1em}%
\begin{subfigure}{.475\linewidth}
  \centering
  \includegraphics[width=\linewidth]{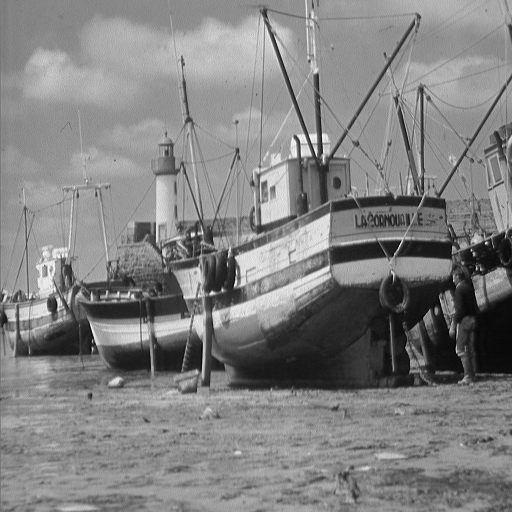}
  \caption{Boat}
  \label{sfig:Boat}
\end{subfigure}
\caption{Images used to compare the proposed algorithm with the BM3D in terms of PSNR, running time, and visual quality}
\label{fig:test.images}
\end{figure}
\begin{figure}[t]
\centering
\begin{subfigure}{.49\linewidth}
  \centering
  \includegraphics[width=\linewidth]{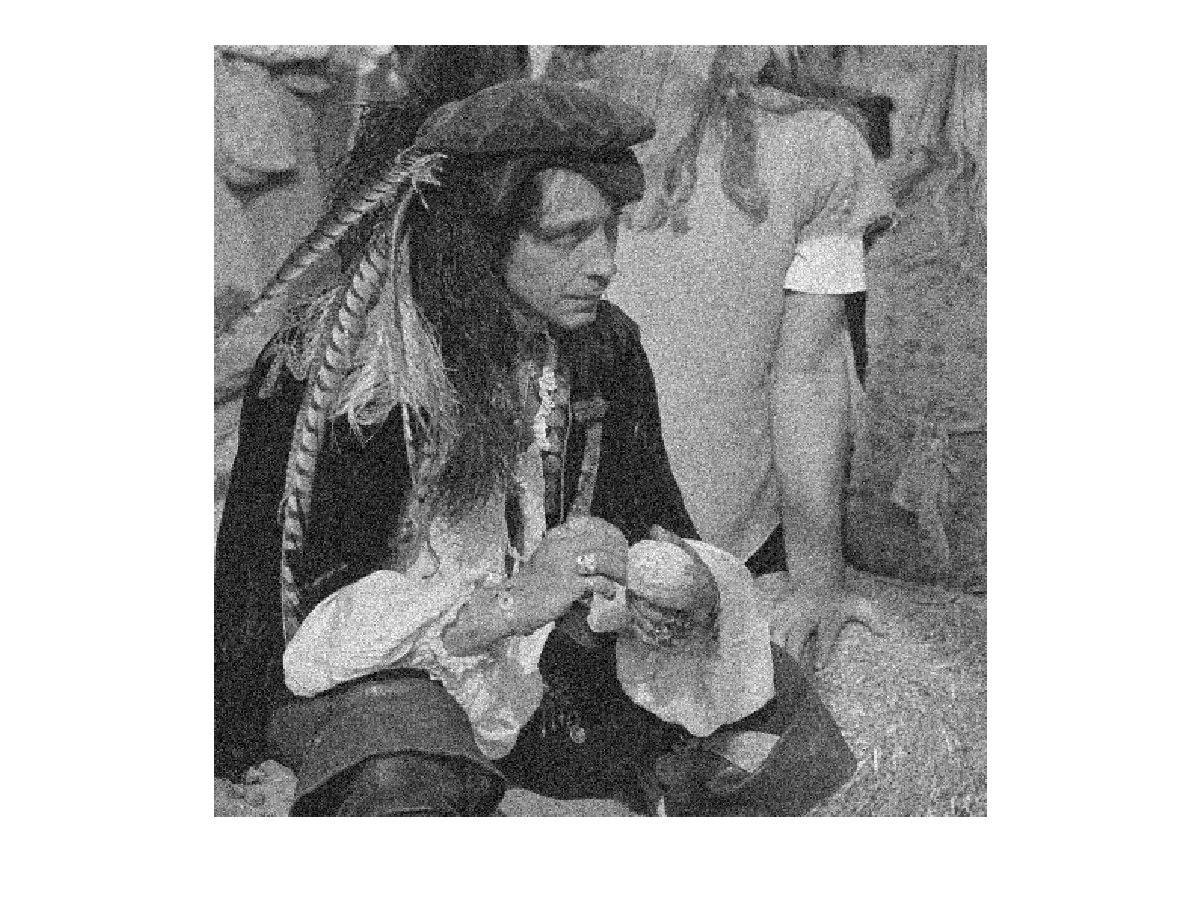}
  \caption{Man: Noisy with $\sigma=20$}
  \label{sfig:Man_std20}
\end{subfigure}%
\hspace{0.5em}%
\begin{subfigure}{.49\linewidth}
  \centering
  \includegraphics[width=\linewidth]{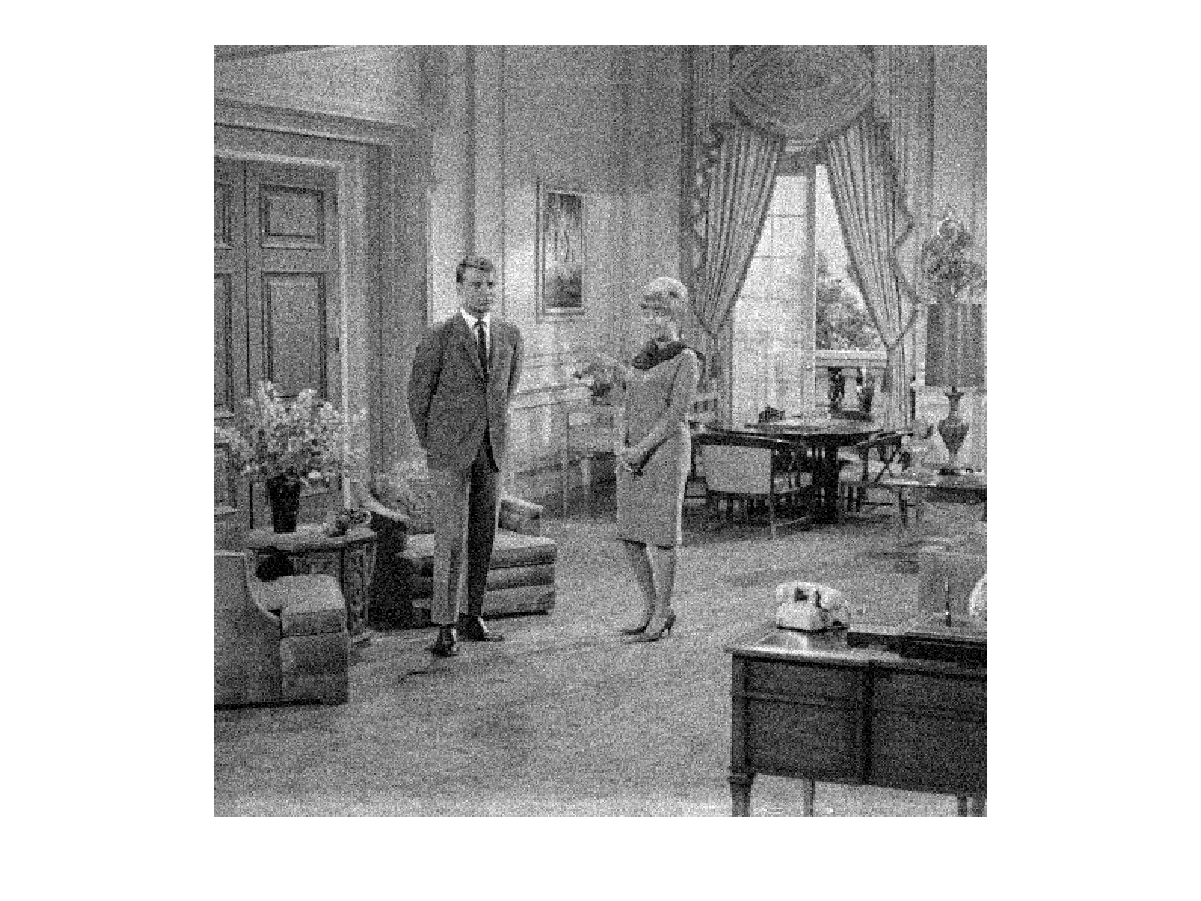}
  \caption{Couple: Noisy with $\sigma=20$}
  \label{sfig:Couple_std20}
\end{subfigure}\\
\vspace{0.25em}
\begin{subfigure}{0.49\linewidth}
  \centering
  \includegraphics[width=\linewidth]{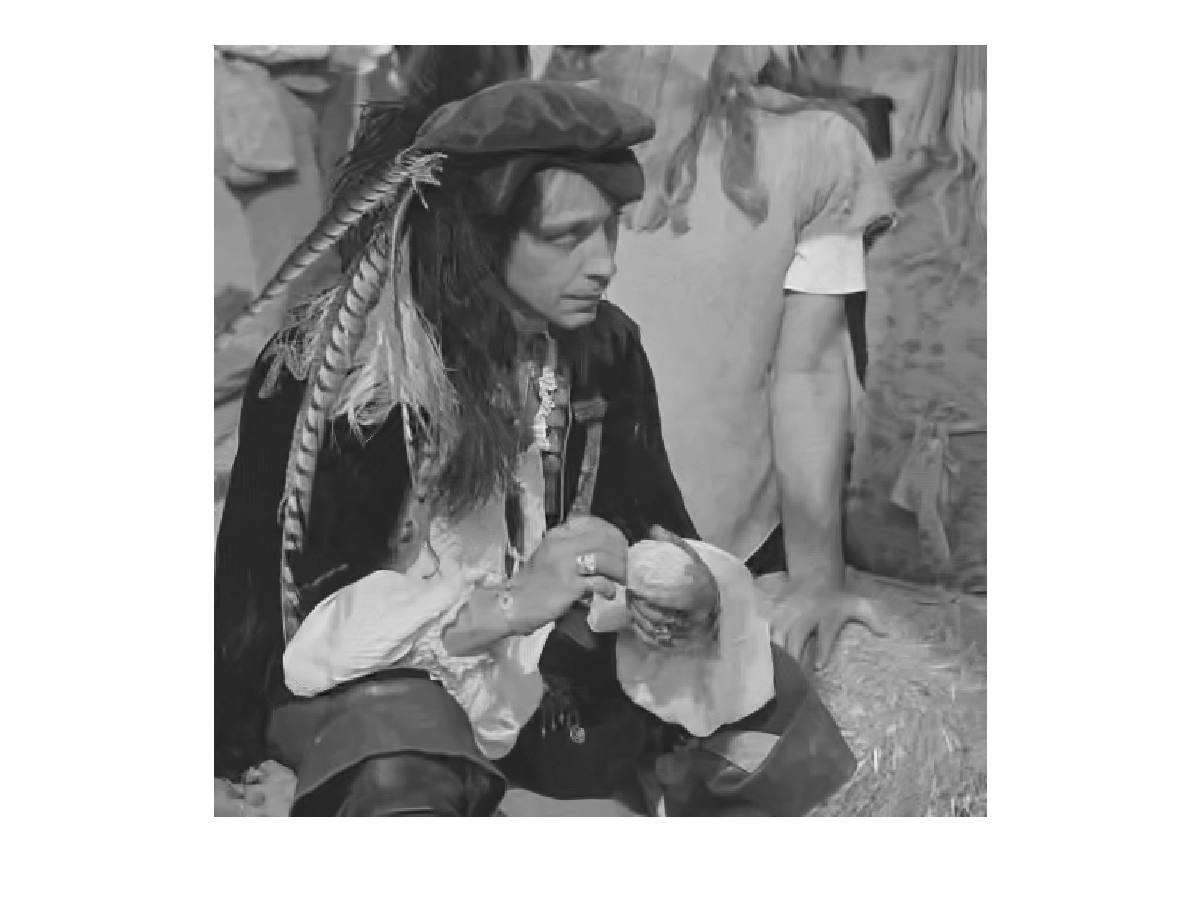}
  \caption{Man: BM3D with PSNR 30.59 dB}
  \label{sfig:Man_std20BM3D}
\end{subfigure}%
\hspace{0.5em}%
\begin{subfigure}{0.49\linewidth}
  \centering
  \includegraphics[width=\linewidth]{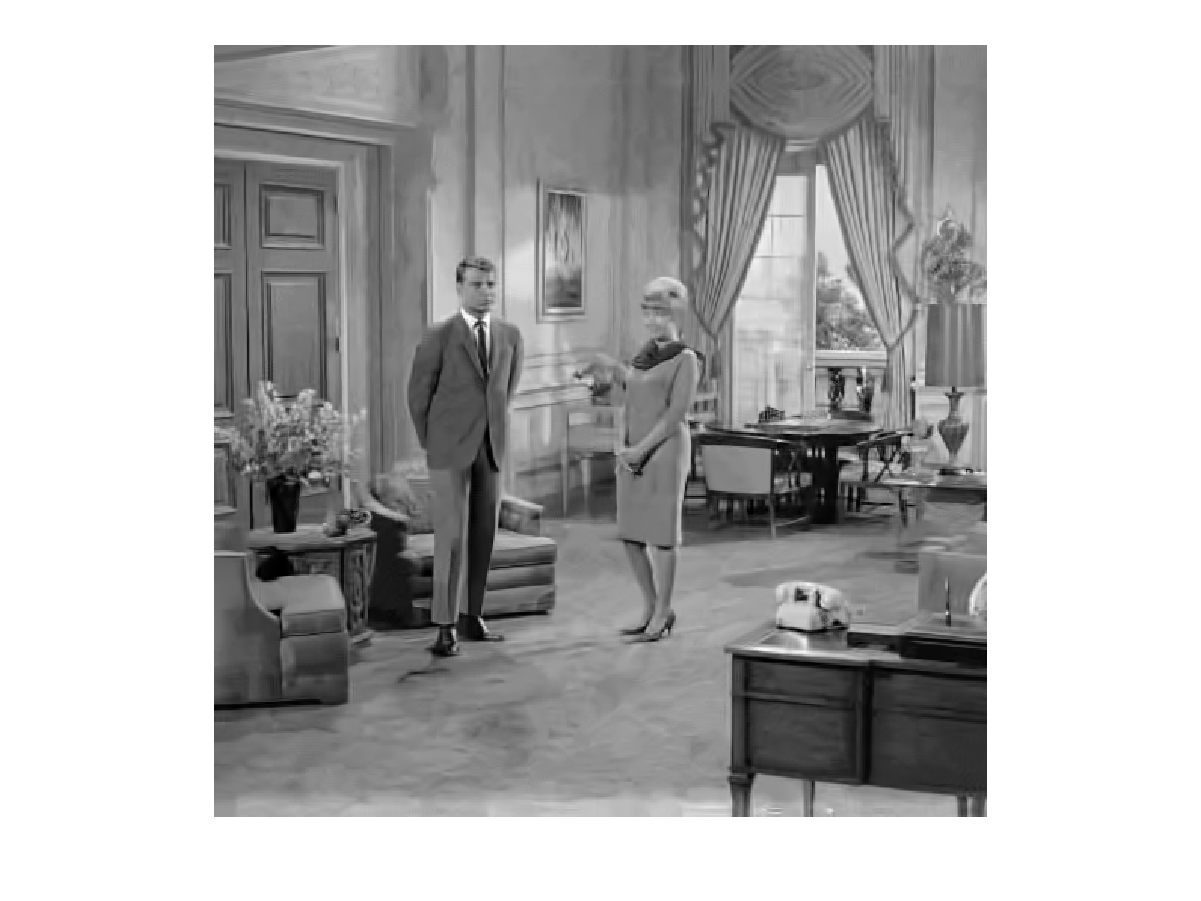}
  \caption{Couple: BM3D with PSNR 30.76 dB}
  \label{sfig:Couple_std20BM3D}
\end{subfigure}\\
\vspace{0.25em}
\begin{subfigure}{0.49\linewidth}
  \centering
  \includegraphics[width=\linewidth]{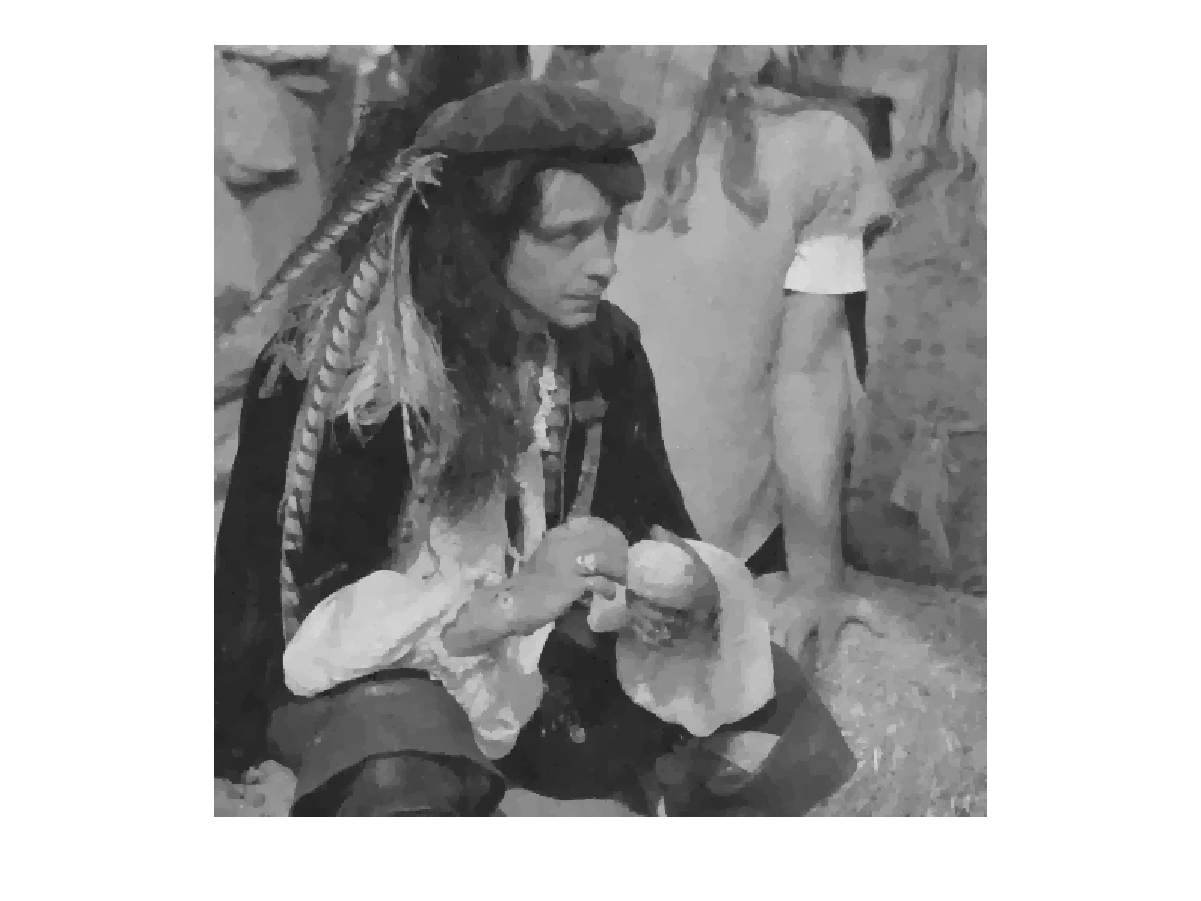}
  \caption{Man: Proposed with PSNR 29.42 dB}
  \label{sfig:Man_std20graph}
\end{subfigure}%
\hspace{0.5em}%
\begin{subfigure}{0.49\linewidth}
  \centering
  \includegraphics[width=\linewidth]{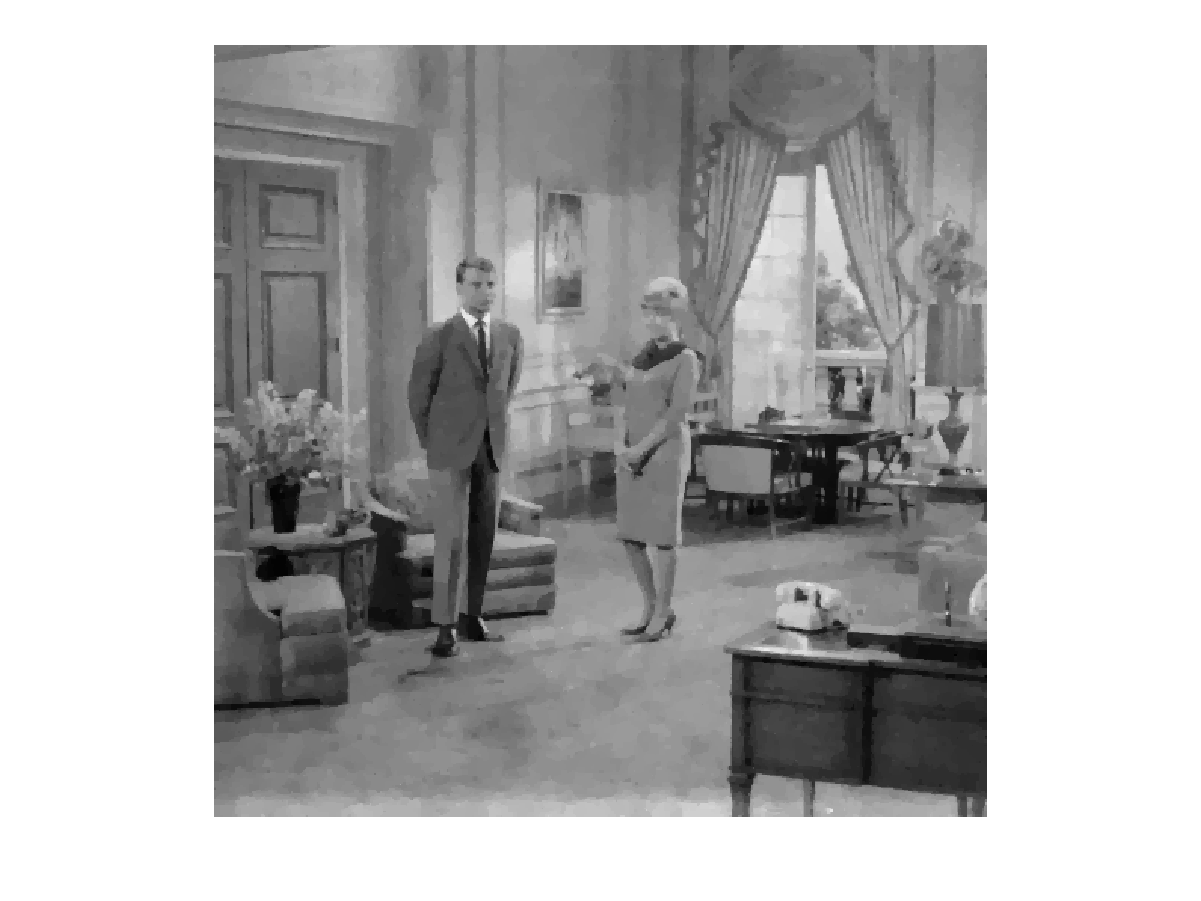}
  \caption{Couple: Proposed with PSNR 28.90 dB}
  \label{sfig:Couple_std20graph}
\end{subfigure}
\caption{Subjective visual quality comparison between denoised images by BM3D and the proposed method. Noise with $\sigma=20$ has been added to the images shown in \ref{sfig:Man} and \ref{sfig:Couple}. In both cases it is clear that though there is a gap in PSNR values, the methods generate reconstructions of similar visual quality. The same conclusion holds for the other images shown in Figure \ref{fig:test.images} and different noise levels.}
\label{fig:test.images.visual}
\end{figure}
In this subsection we compare numerically the proposed algorithm with the BM3D denoising method introduced in \cite{BM3D1}. This method is a state-of-the-art filtering method applicable to denoising among other important problems in image processing. It is therefore interesting to make a comparison although BM3D is not directly designed to solve the TV minimization problem.
\par In our experiments we used different types of images and Gaussian noise with standard deviation ranging from 5 to 100. More precisely, we have reproduced and compared results obtained by BM3D and by the proposed method by testing all images available at the online repository \url{http://www.cs.tut.fi/~foi/GCF-BM3D/}. However, since all results obtained lead to the same conclusions, we have only reported results obtained for test images illustrated in Figure \ref{fig:test.images}. Comparison with BM3D has been made in terms of PSNR, running time and visual quality.
\begin{itemize}
\item The implementation of the BM3D was obtained from \cite{BM3D1} and the online repository \url{http://www.cs.tut.fi/~foi/GCF-BM3D/}.
\vskip1mm  
\item We implemented the proposed Algorithm \ref{ALGO} and optimized the regularization parameter $t$ by performing a brute-force optimization with uniformly quantized values in suitable sub-intervals of the interval $[1,100]$ corresponding to noise levels ranging from 5, 10, $\ldots$, 100. 
\end{itemize}
Obtained results in terms of PSNR, running time and visual quality are shown in Table \ref{PSNR:table}, Table \ref{Time:table} and Figure \ref{fig:test.images.visual} respectively. As might be expected, the PSNR values for the test images is higher for the BM3D method as shown in Table \ref{PSNR:table}. The difference in PSNR values is in the range of 1 to 2 dB. In terms of running time, as can be seen in Table \ref{Time:table}, the proposed algorithm is the most competitive in the low noise regime. However, it is also more sensitive to the noise level than the BM3D method. In terms of visual quality, the methods produce comparable results.
\begin{table}[h]
\centering
\resizebox{\columnwidth}{!}{%
\begin{tabular}{|c|c|c|c|c|c|c|c|c|} 
\hline
\textbf{} & \multicolumn{2}{c|}{\textbf{\begin{tabular}[c]{@{}c@{}}Man\\ $512\times 512$\end{tabular}}} & \multicolumn{2}{c|}{\textbf{\begin{tabular}[c]{@{}c@{}}Couple\\ $512\times 512$\end{tabular}}} & \multicolumn{2}{c|}{\textbf{\begin{tabular}[c]{@{}c@{}}Hill\\ $512\times 512$\end{tabular}}} & \multicolumn{2}{c|}{\textbf{\begin{tabular}[c]{@{}c@{}}Boat\\ $512\times 512$\end{tabular}}} \\ \hline
\textbf{$\sigma$} & \textbf{BM3D} & \textbf{Proposed} & \textbf{BM3D} & \textbf{Proposed} & \textbf{BM3D} & \textbf{Proposed} & \textbf{BM3D} & \textbf{Proposed} \\ \hline
\textbf{5}   & 37.82 & 36.80 & 37.52 & 36.50 & 37.13 & 36.40 & 37.06 & 36.27  \\ \hline
\textbf{10}  & 33.98 & 32.77 & 34.04 & 32.44 & 33.62 & 32.57 & 33.89 & 32.49  \\ \hline
\textbf{15}  & 31.93 & 30.73 & 32.11 & 30.31 & 31.86 & 30.70 & 32.10 & 30.50  \\ \hline
\textbf{20}  & 30.59 & 29.34 & 30.76 & 28.90 & 30.72 & 29.50 & 30.83 & 29.17  \\ \hline
\textbf{25}  & 29.62 & 28.48 & 29.72 & 27.90 & 29.85 & 28.63 & 29.85 & 28.18  \\ \hline
\textbf{30}  & 28.86 & 27.75 & 28.87 & 27.09 & 29.15 & 27.94 & 29.04 & 27.41  \\ \hline
\textbf{35}  & 28.22 & 27.16 & 28.15 & 26.45 & 28.49 & 27.39 & 28.32 & 26.78  \\ \hline
\textbf{40}  & 27.65 & 26.67 & 27.48 & 25.93 & 27.88 & 26.94 & 27.63 & 26.25  \\ \hline
\textbf{45}  & 27.17 & 26.25 & 26.91 & 25.09 & 27.43 & 26.56 & 27.12 & 25.81  \\ \hline
\textbf{50}  & 26.81 & 25.88 & 26.46 & 24.75 & 27.19 & 26.23 & 26.67 & 25.43  \\ \hline
\textbf{55}  & 26.44 & 25.56 & 26.01 & 24.48 & 26.74 & 25.94 & 26.27 & 25.10  \\ \hline
\textbf{60}  & 26.14 & 25.27 & 25.66 & 24.45 & 26.52 & 25.68 & 25.90 & 24.80  \\ \hline
\textbf{65}  & 25.90 & 25.01 & 25.29 & 24.18 & 26.12 & 25.44 & 25.56 & 24.52  \\ \hline
\textbf{70}  & 25.56 & 24.77 & 25.00 & 23.94 & 25.93 & 25.23 & 25.25 & 24.27  \\ \hline
\textbf{75}  & 25.32 & 24.56 & 24.70 & 23.72 & 25.68 & 25.03 & 24.97 & 24.04  \\ \hline
\textbf{80}  & 25.06 & 24.36 & 24.42 & 23.52 & 25.43 & 24.85 & 24.70 & 23.83  \\ \hline
\textbf{85}  & 24.86 & 24.17 & 24.21 & 23.33 & 25.10 & 24.68 & 24.45 & 23.64  \\ \hline
\textbf{90}  & 24.63 & 24.00 & 23.94 & 23.16 & 24.98 & 24.52 & 24.22 & 23.46  \\ \hline
\textbf{95}  & 24.39 & 23.83 & 23.67 & 23.01 & 24.66 & 24.38 & 24.01 & 23.29  \\ \hline
\textbf{100} & 24.22 & 23.67 & 23.51 & 22.86 & 24.58 & 24.24 & 23.80 & 23.13  \\ \hline
\end{tabular}%
}
\caption{PSNR (dB) results of the proposed method and BM3D method}
\label{PSNR:table}
\end{table}
\begin{table}[h]
\centering
\resizebox{\columnwidth}{!}{%
\begin{tabular}{|c|c|c|c|c|c|c|c|c|}
\hline
\textbf{} & \multicolumn{2}{c|}{\textbf{\begin{tabular}[c]{@{}c@{}}Man\\ $512\times 512$\end{tabular}}} & \multicolumn{2}{c|}{\textbf{\begin{tabular}[c]{@{}c@{}}Couple\\ $512\times 512$\end{tabular}}} & \multicolumn{2}{c|}{\textbf{\begin{tabular}[c]{@{}c@{}}Hill\\ $512\times 512$\end{tabular}}} & \multicolumn{2}{c|}{\textbf{\begin{tabular}[c]{@{}c@{}}Boat\\ $512\times 512$\end{tabular}}} \\ \hline
\textbf{$\sigma$} & \textbf{BM3D} & \textbf{Proposed} & \textbf{BM3D} & \textbf{Proposed} & \textbf{BM3D} & \textbf{Proposed} & \textbf{BM3D} & \textbf{Proposed} \\ \hline
\textbf{5}   & 7.50  & 1.24  & 7.40  & 0.99   & 8.10  & 0.98   &  11.4   & 0.99    \\ \hline
\textbf{10}  & 8.00  & 2.37  & 7.80  & 2.04   & 8.50  & 2.15   &  8.4    & 2.21    \\ \hline
\textbf{15}  & 8.40  & 3.40  & 8.20  & 3.09   & 8.90  & 3.52   &  8.7    & 3.28    \\ \hline
\textbf{20}  & 8.60  & 4.66  & 8.50  & 4.06   & 9.20  & 4.67   &  11.6   & 4.15    \\ \hline
\textbf{25}  & 8.90  & 5.42  & 8.70  & 5.14   & 9.50  & 5.72   &  9.0    & 5.08    \\ \hline
\textbf{30}  & 9.00  & 6.47  & 9.00  & 6.05   & 9.50  & 6.89   &  9.1    & 6.23    \\ \hline
\textbf{35}  & 9.10  & 7.57  & 9.00  & 7.02   & 9.50  & 7.77   &  9.0    & 6.94    \\ \hline
\textbf{40}  & 8.80  & 8.55  & 8.60  & 7.99   & 9.30  & 8.96   &  8.8    & 8.18    \\ \hline
\textbf{45}  & 10.50 & 9.69  & 10.50 & 8.88   & 10.70 & 9.89   &  10.6   & 9.01    \\ \hline
\textbf{50}  & 10.60 & 10.84 & 10.60 & 9.66   & 10.80 & 10.70  &  10.7   & 9.62    \\ \hline
\textbf{55}  & 10.70 & 11.62 & 10.70 & 10.42  & 10.90 & 12.12  &  10.8   & 10.54   \\ \hline
\textbf{60}  & 10.80 & 12.46 & 10.70 & 11.48  & 11.00 & 13.20  &  10.9   & 11.54   \\ \hline
\textbf{65}  & 10.80 & 14.00 & 10.80 & 12.14  & 11.10 & 13.61  &  11.0   & 12.42   \\ \hline
\textbf{70}  & 10.90 & 14.35 & 10.90 & 13.26  & 11.10 & 15.13  &  13.4   & 13.15   \\ \hline
\textbf{75}  & 11.00 & 15.53 & 11.00 & 13.87  & 11.20 & 15.46  &  11.3   & 14.30   \\ \hline
\textbf{80}  & 11.00 & 16.60 & 11.10 & 15.08  & 11.20 & 16.02  &  11.3   & 15.20   \\ \hline
\textbf{85}  & 11.10 & 17.40 & 11.10 & 16.02  & 11.40 & 17.65  &  11.4   & 15.85   \\ \hline
\textbf{90}  & 11.20 & 17.99 & 11.20 & 16.76  & 11.40 & 18.34  &  11.4   & 16.44   \\ \hline
\textbf{95}  & 11.20 & 18.97 & 11.30 & 17.93  & 11.50 & 19.54  &  11.4   & 17.22   \\ \hline
\textbf{100} & 11.20 & 19.32 & 11.30 & 18.32  & 11.50 & 20.49  &  11.4   & 17.38   \\ \hline
\end{tabular}%
}
\caption{Running time (seconds) results of the proposed method and BM3D method}
\label{Time:table}
\end{table}
\section{Conclusion}\label{SecConclznz}
\par In this work we proposed an iterative algorithm for total variation minimization on graphs and proved its convergence. The algorithm that is presented can be viewed as a coordinate descent on dual space and can be run on a parallel computer architecture, which makes it suitable to handle large graphs and data sets. The algorithm is simple, easy to implement and converges to the exact minimizer with fewer iterations compared to the Split-Bregman and Primal-Dual algorithms. Furthermore, in order to compare the proposed method with other state-of-the-art denoising methods, BM3D was chosen as a benchmark and obtained results still show competitive performance for the proposed method. In a follow-up work we intend to further study the convergence rate of the algorithm and include additional imaging scenarios.
\section*{Acknowledgements}
The authors are grateful to the anonymous referees for insightful comments that helped to improve the paper. The research of Japhet Niyobuhungiro was partially supported by The World Academy of Sciences (TWAS), RGA No. 17-356\ RG/MATHS/AF/AC\textunderscore I - FR3240297744. Eric Setterqvist acknowledges support by the Austrian Science Fund (FWF) within the national research network S117 ``Geometry $+$ Simulation", subproject 4.
\bibliographystyle{plain}
\bibliography{Mainbibfile}

\end{document}